\documentclass[a4paper,english,11pt,oneside]{article}
\usepackage[a4paper]{geometry}
\usepackage{ifthen}
\usepackage[english]{babel}
\usepackage{bibgerm}       		
\usepackage[utf8]{inputenc}
\usepackage[hyphens]{url}           		
\usepackage{hyperref} 				
\usepackage{caption}
\usepackage{subcaption}
\usepackage{amsmath}
\usepackage{amsfonts}
\usepackage{amssymb}
\usepackage{mathtools}
\usepackage{siunitx}
 \DeclareSIUnit\bar{bar}
\usepackage[linesnumbered, ruled, noend]{algorithm2e}
\usepackage{graphicx}
\usepackage{enumerate}
\usepackage{xcolor}
\usepackage{xspace}

\usepackage{times}

\usepackage{tikz}
\usetikzlibrary{patterns}
\usetikzlibrary{positioning}

\usetikzlibrary{shapes}
\usetikzlibrary{shapes.geometric}
\usetikzlibrary{backgrounds}
\usetikzlibrary{arrows.meta}

\usepackage{footnote}
\makesavenoteenv{tabular}
\makesavenoteenv{table}
\usepackage{threeparttable}
\usepackage{pgfplots}

\providecommand{\keywords}[1]
{
  \small	
  \textbf{\textit{Keywords---}} #1
}

\newcommand{\graph}{\ensuremath{G}\xspace}
\newcommand{\generalGraph}{\ensuremath{\graph=(\setVertices,\setArcs)}\xspace}

\newcommand{\temperature}{\ensuremath{T}\xspace}
\newcommand{\logTemp}{\ensuremath{\Delta\temperature_\mathrm{ln}}\xspace}
\newcommand{\outTemp}{\ensuremath{\temperature^{\mathrm{out}}}\xspace}
\newcommand{\outTempArc}[1]{\ensuremath{\temperature_{#1}^{\mathrm{out}}}\xspace}
\newcommand{\soilTemp}{\ensuremath{\temperature^{\mathrm{s}}}\xspace}
\newcommand{\inTemp}{\ensuremath{\temperature^{\mathrm{in}}}\xspace}
\newcommand{\inTempArc}[1]{\ensuremath{\temperature_{#1}^{\mathrm{in}}}\xspace}
\newcommand{\pressure}{\ensuremath{p}\xspace}
\newcommand{\outPressure}{\ensuremath{\pressure^{\mathrm{out}}}\xspace}
\newcommand{\inPressure}{\ensuremath{\pressure^{\mathrm{in}}}\xspace}
\newcommand{\heatFlow}{\ensuremath{Q}\xspace}
\newcommand{\enthalpy}{\ensuremath{H}\xspace}
\newcommand{\massFlow}{\ensuremath{q}\xspace}
\newcommand{\diameter}{\ensuremath{d}\xspace}
\newcommand{\length}{\ensuremath{L}\xspace}
\newcommand{\outerDiameter}{\ensuremath{\diameter_{\mathrm{o}}}\xspace}
\newcommand{\pipeDepth}{\ensuremath{s}\xspace}
\newcommand{\pipeLossCoeff}{\ensuremath{\phi}\xspace}
\newcommand{\density}{\ensuremath{\rho}\xspace}
\newcommand{\dynViscosity}{\ensuremath{\nu}\xspace}
\newcommand{\frictionFactor}{\lambda}
\newcommand{\heatCap}{\ensuremath{c_\mathrm{p}}\xspace}
\newcommand{\heatCapArc}[1]{\ensuremath{c_{\mathrm{p},#1}}\xspace}
\newcommand{\thermCond}{\ensuremath{\kappa}\xspace}
\newcommand{\heatTransCoeff}{\ensuremath{\alpha}\xspace}
\newcommand{\heatTransRate}{\ensuremath{k}\xspace}
\newcommand{\jouleThomson}{\ensuremath{\mu_{\mathrm{JT}}}\xspace}
\newcommand{\cost}{\ensuremath{c}\xspace}
\newcommand{\pressureTol}{\ensuremath{\epsilon_\mathrm{\pressure}}\xspace}

\newcommand{\temperatureTol}{\ensuremath{\epsilon_\mathrm{\temperature}}\xspace}
\newcommand{\frictionTol}{\ensuremath{\epsilon_\mathrm{\frictionFactor}}\xspace}

\newcommand{\setVertices}{\ensuremath{V}\xspace}
\newcommand{\setEntries}{\ensuremath{\setVertices^{+}}\xspace}
\newcommand{\setExits}{\ensuremath{\setVertices^{-}}\xspace}
\newcommand{\setBoundaryNodes}{\ensuremath{\setVertices^\mathrm{b}}\xspace}
\newcommand{\setBoundaryNodesAll}{\ensuremath{\setBoundaryNodes = \setEntries \dot\cup\,\setExits}\xspace} 
\newcommand{\setInnerNodes}{\ensuremath{\setVertices^{0}}\xspace}
\newcommand{\setMiddleNodes}{\ensuremath{\setVertices^{m}}\xspace}
\newcommand{\setLeafs}{\ensuremath{L}\xspace}
\newcommand{\setVerticesAll}{\ensuremath{\setVertices = \setEntries \dot\cup\,\setExits \dot\cup\, \setInnerNodes}\xspace}

\newcommand{\vertex}{\ensuremath{v}\xspace}
\newcommand{\anotherVertex}{\ensuremath{u}\xspace}
\newcommand{\vertexNumbered}[1]{\ensuremath{\vertex_{#1}}\xspace}

\newcommand{\vertexPressure}[1]{\ensuremath{\pressure_{#1}}\xspace}
\newcommand{\boundaryVal}[1]{\ensuremath{b_{#1}}\xspace}

\newcommand{\setArcs}{\ensuremath{A}\xspace}
\newcommand{\arc}{\ensuremath{a}\xspace}

\newcommand{\arcCost}[2]{\ensuremath{\cost_{#1,#2}}\xspace}

\newcommand{\arcFlow}[1]{\ensuremath{\massFlow_{#1}}\xspace}
\newcommand{\pipeDiameter}[1]{\ensuremath{\diameter_{#1}}\xspace}
\newcommand{\setPipes}{\ensuremath{\setArcs^{\mathrm{pi}}}\xspace}
\newcommand{\setDiameter}[1]{\ensuremath{D_{#1}}\xspace}
\newcommand{\setPumps}{\ensuremath{\setArcs^{\mathrm{pu}}}\xspace}

\newcommand{\setIncomingArcs}[1]{\ensuremath{\delta^{-}(#1)}\xspace}
\newcommand{\setOutgoingArcs}[1]{\ensuremath{\delta^{+}(#1)}\xspace}
\newcommand{\setIncomingArcsDef}[2]{\ensuremath{\delta^{-}(#1) = \lbrace \arc \in \setArcs : a = (#2, #1) \rbrace}\xspace}
\newcommand{\setOutgoingArcsDef}[2]{\ensuremath{\delta^{+}(#1) = \lbrace \arc \in \setArcs : a = (#1, #2) \rbrace}\xspace}

\newcommand{\Rpos}{\ensuremath{\mathbb{R}_{\ge0}}\xspace}
\newcommand{\UB}[2]{\ensuremath{\overline{#1}_{#2}}\xspace}
\newcommand{\LB}[2]{\ensuremath{\underbar{$#1$}_{#2}}\xspace}

\newcommand{\arcLength}{\ensuremath{\length_\arc}\xspace}

\newcommand{\coTwo}{\ensuremath{\mathrm{CO}_2}\xspace}
\newcommand\Rey{\ensuremath{\mathrm{Re}}\xspace}
\newcommand\Pra{\ensuremath{\mathrm{Pr}}\xspace}
\newcommand\Nus{\ensuremath{\mathrm{Nu}}\xspace}

\newcommand{\python}{{\sc python}\xspace}
\newcommand{\gurobi}{{\sc gurobi}\xspace}
\newcommand{\gurobiVersion}[2]{{\sc gurobi}~\oldstylenums{#1.#2}\xspace}

\newcommand{\gascalcVersion}[1]{{\sc gascalc}~\oldstylenums{#1}\xspace}

\DeclareRobustCommand{\markerentry}{\raisebox{0.5pt}{\tikz{\node[draw,scale=0.5,regular polygon,regular polygon sides=6,fill=green, rotate=90](){};}}}
\DeclareRobustCommand{\markerexit}{\raisebox{0.5pt}{\tikz{\node[draw,scale=0.5,regular polygon,regular polygon sides=6,fill=red](){};}}}
\DeclareRobustCommand{\markerpump}{\raisebox{0pt}{\tikz{\node[draw,scale=0.5,diamond,fill=orange](){};}}}
\DeclareRobustCommand{\markerinnode}{\raisebox{0pt}{\tikz{\node[draw,scale=0.5,circle,fill=blue](){};}}}
\allowdisplaybreaks

\usepackage{csquotes}
\usepackage{zibtitlepage}
\usepackage{authblk}

\begin{document}


\title{Optimal discrete pipe sizing for tree-shaped $\text{CO}_2$ networks}
%
\author[1]{Jaap Pedersen\footnote{corresponding author, pedersen@zib.de, \ZTPOrcid{0000-0003-4047-0042}}}
\author[1,2]{Thi Thai Le\footnote{\ZTPOrcid{0000-0001-7886-4878}}}
\author[1,2]{Thorsten Koch\footnote{\ZTPOrcid{0000-0002-1967-0077}}}
\author[1]{Janina Zittel\footnote{\ZTPOrcid{0000-0002-0731-0314}}}

\affil[1]{Zuse Institute Berlin, Berlin, Germany}
\affil[2]{Technische Universit{\"a}t Berlin, Berlin, Germany}

\maketitle

\begin{abstract}
    Many energy-intensive industries, like the steel industry, plan to switch to renewable energy sources. Other industries, such as the cement industry, have to rely on carbon capture utilization and storage (CCUS) technologies to reduce their production processes' inevitable carbon dioxide (\coTwo) emissions. However, a new transport infrastructure needs to be established to connect the point of capture and the point of storage or utilization. Given a tree-shaped network transporting captured \coTwo from multiple sources to a single sink, we investigate how to select optimal pipeline diameters from a discrete set of diameters. The general problem of optimizing arc capacities in potential-based fluid networks is already a challenging mixed-integer nonlinear optimization problem. The problem becomes even more complex when adding the highly sensitive nonlinear behavior of \coTwo regarding temperature and pressure changes. We propose an iterative algorithm that splits the problem into two parts: a) the pipe-sizing problem under a fixed supply scenario and temperature distribution and b) the thermophysical modeling, including mixing effects, the Joule-Thomson effect, and the heat exchange with the surrounding environment. We show the effectiveness of our approach by applying our algorithm to a real-world network planning problem for a \coTwo network in Germany.
\end{abstract}


\section{Introduction}\label{sec:introduction}
Industry contributes between \SI{20}-\SI{30}{\%} to the global greenhouse gas emissions \cite{IPCC2021, iea2020a}. Whereas energy-intensive industries like steel production may reduce their carbon footprint by switching towards renewable energy sources, it is harder to decarbonize for other industries like the cement production where the carbon dioxide (\coTwo) emission arise from the chemical process itself \cite{rissman2020, davis2018}. However, to fulfill the Paris Agreement with a global average temperature increase ``well below \SI{2}{\celsius} and towards \SI{1.5}{\celsius}'', all sectors have to reach net-zero by 2050 \cite{paris2015}. Besides an increasing material efficiency and an end-to-end life cycle, carbon capture and storage (CCS) and utilization (CCU) technologies are becoming more important \cite{davis2018, iea2019, iea2020, bui2018, harvey2021}. For the cost-efficient deployment of CCS and CCU, \coTwo transport is the linking component between the point of capture and the point of storage and utilization. Therein, pipelines and pipeline networks have proven to be the cheapest option for transporting \coTwo in large quantities onshore and, depending on the distance, and the amount, offshore, \cite{iea2020, JRC2011, ipcc2005}. This is also in line with other supply chain models, e.g., in the context of hydrogen transport networks \cite{reuss2019}.

In this study, a tree-shaped network is investigated for transporting captured \coTwo. The goal is to find the optimal pipeline diameters from a discrete set of diameters for a given network topology, i.e., pipe lengths, pipe route, and location of the pumping stations are known. Optimizing arc capacities in a network is a common problem in many real-world applications such as water, gas, hydrogen, or in our case, \coTwo transportation networks. The problem is an extension of the potential-based network flow problem, in which the relationship between arc flows and potential, e.g., gas pressure or water height, is described by physical laws. These relationships are usually nonlinear. For an overview of how to model and solve the potential-based network flow problem, we refer to \cite{gross2019} and the references therein. In combination with binary control decisions of the active elements in the network, e.g., closing or opening of a valve in gas networks, finding the optimal solution results in challenging mixed-integer nonlinear optimization problems (MINLP). 

In \cite{bragalli2012} and \cite{dambrosio2015}, the problem of finding optimal pipeline diameters in water networks is solved by applying mixed-integer nonlinear programming techniques. The authors of \cite{dambrosio2015} also consider active elements, such as pumps and valves. These references deal with networks with any kind of topology, i.e., meshed or tree-shaped. In \cite{shiono2016}, optimal pipe-sizing in tree-shaped gas networks with a single source and given demands is discussed. Here, algorithms are proposed for approximating discrete diameters of pipelines based on previously found continuous diameters. The process of expanding existing gas networks by building new parallel pipelines, also known as \emph{looping}, is studied, for example, in \cite{andre2009} and \cite{lenz2022}. The authors of \cite{andre2009} propose a two-stage approach by firstly identifying the critical pipelines and secondly deciding on the discrete diameters using Branch-and-Bound. Besides a novel approach for continuous looping, \cite{lenz2022} also proposes a method to solve the problem for a multi-case scenario which helps practitioners prepare for different future demand and flow scenarios. 

The case of multi-case scenarios is also discussed in \cite{robinius2019}, in which a discrete pipe-sizing problem is studied in the context of hydrogen networks. Here, the authors consider a tree-shaped single-source network with known demand. By creating a finite representation of an infinite set of scenarios, they are able to solve the problem as a standard mixed-integer linear optimization problem (MILP).

\begin{figure}
    \centering
    \begin{subfigure}[t]{.45\linewidth}
        \includegraphics[width=\linewidth]{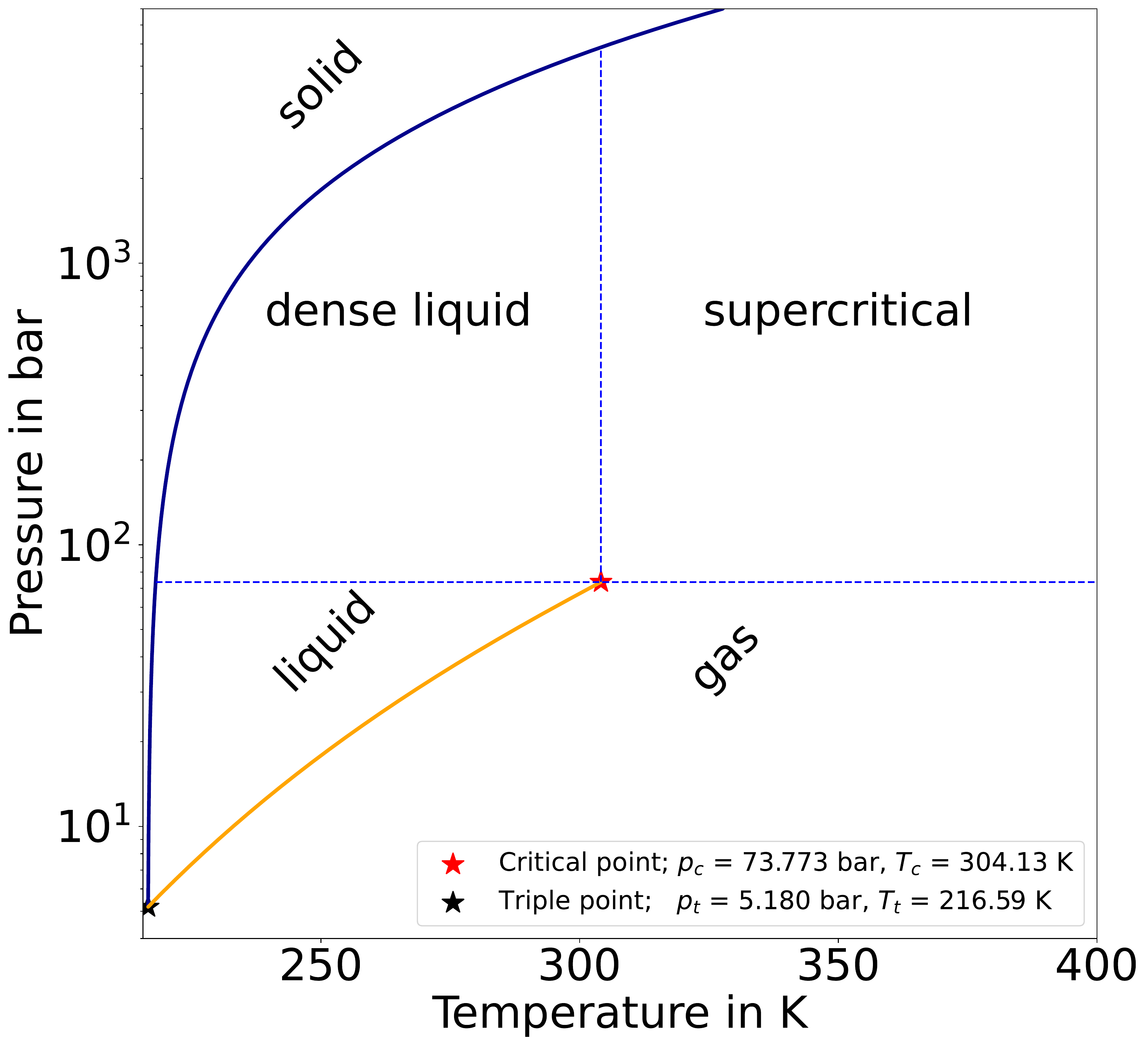}
    \caption{\coTwo phase diagram, data from \cite{bell2014}}
    \label{fig:phaseDiagram}
    \end{subfigure}
    \begin{subfigure}[t]{.5\linewidth}
        \includegraphics[width=\linewidth]{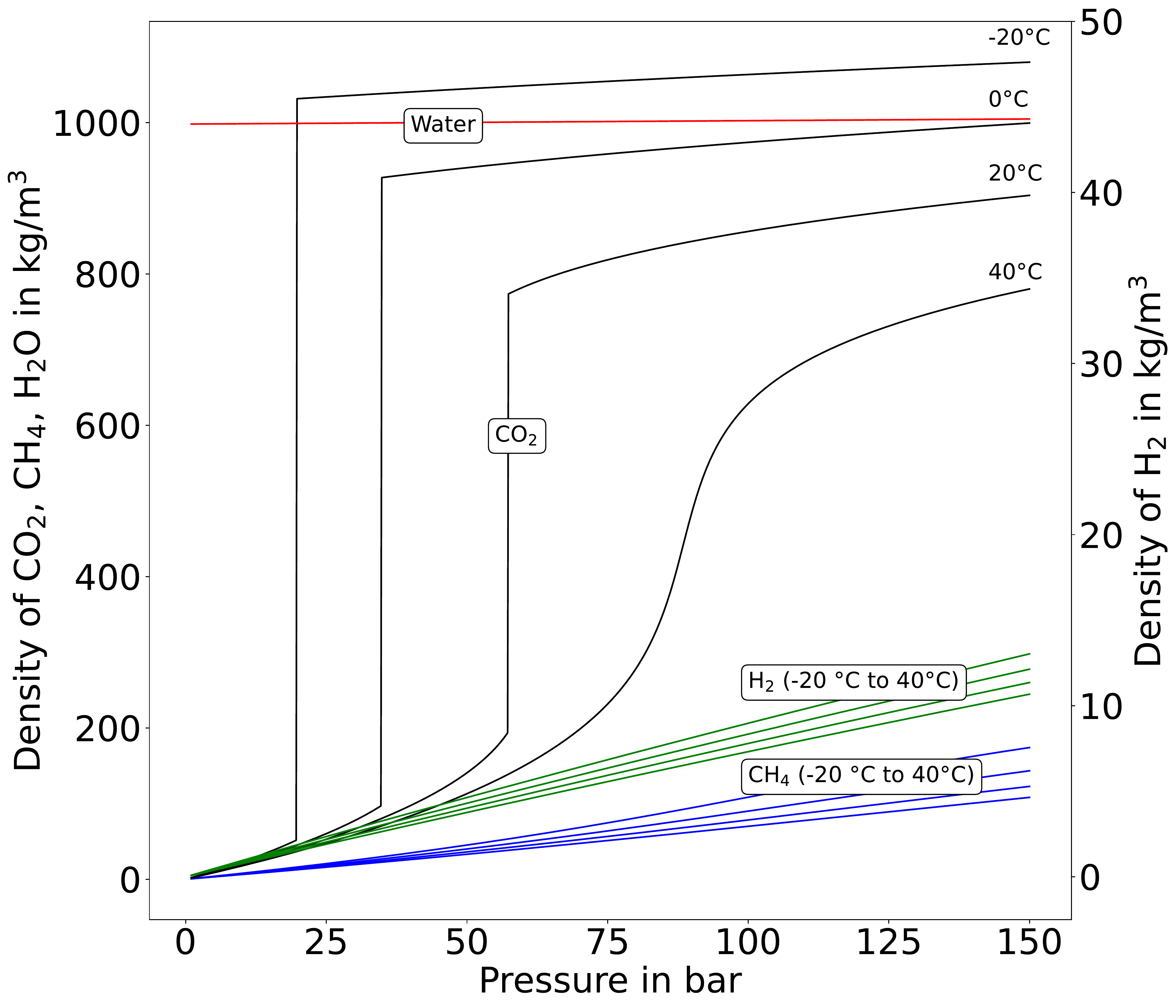}
    \caption{Density of different fluids,  similar to \cite{nimtz2016}, data from \cite{bell2014}}
    \label{fig:compare_co2}
    \end{subfigure}
    \caption{}
    \label{fig:co2_data}
\end{figure}

In the references mentioned above, the authors investigate networks in which no phase changes of the transported fluid occur, allowing them to make valid assumptions regarding the physics of the system. For example, one can easily assume that liquid water is an incompressible fluid, i.e., density does not change with temperature and pressure. Also, it is fair to say that gas transport is isothermic and that in the range of the considered pressure levels, the physical properties, e.g., the dynamic viscosity, are constant or a linear function in terms of pressure. However, \coTwo is usually transported in liquified or supercritical\footnote{In supercritical form, the liquid and gaseous phase of the fluid cannot be distinct, as both have the same density} form. In addition to a possible transition from supercritical to liquid, the physical properties of \coTwo are also highly sensitive to pressure and temperature changes within the liquid phase. Furthermore, the relationships between these properties and pressure and temperature are nonlinear. Figure~\ref{fig:phaseDiagram} shows the \pressure-\temperature phase diagram for \coTwo.
As an example of the more complex behavior of \coTwo compared to other transported fluids mentioned here, Figure~\ref{fig:compare_co2} shows the density for water, \coTwo, hydrogen, and methane, which is the primary content of natural gas, in the relevant pressure ranges. For other quantities, the reader is referred to \cite{nimtz2016}.

For a general understanding of the behavior and physical properties of \coTwo and its transport, we refer to the survey papers \cite{JRC2011, munkejord2016} and the references therein. In the review paper \cite{peletiri2018}, the authors focus on the process of designing a \coTwo pipeline. The authors include pipeline route selection, fluid flow rates, and pipeline diameter and pressure drop calculations. There exist several studies analyzing the simulation of \coTwo transport in a single pipeline. For example, \cite{aursand2013} gives an overview of the state-of-the-art for modeling the transient flow of \coTwo in pipelines. 
In \cite{nimtz2010}, a combined system of a \coTwo pipeline and storage is simulated, showing small changes in boundary conditions will lead to undesirable behavior in the system, e.g., change of phase. A thorough simulation of a complete CCS system for \coTwo-rich fluids is performed in \cite{nimtz2016} consisting of point of capture, transport via a single pipeline, and the injection into an undersea storage. The author also includes simulation of the filling and operation of the pipeline, running-ductile fracture\footnote{A dangerous incident in which an initial small crack may start running alongside the pipeline and, in the worst case, may split the pipeline.}, and changing load profiles. The operation of a multi-source \coTwo pipeline network with fixed diameters and different \coTwo compositions is modeled in \cite{brown2015}. The current literature regarding \coTwo transport relies heavily on the simulation of single pipelines. Optimal pipeline sizing is usually solved by iteratively and manually adapting the choice of diameter with respect to mass flow and allowed pressure drop.  

In this paper, we propose a new method to solve the arc capacity optimization problem taking the challenging thermophysical behavior of \coTwo into account. Due to the complex thermophysical behavior of \coTwo, the problem is distinctly more difficult than the capacity optimization for transporting the other mentioned fluids, i.e., water, natural gas, or hydrogen. With our approach, it is not only possible to find the optimal pipeline diameters for a given inflow scenario but also the pressure and temperature distribution in the network. Furthermore, larger pressure drops are allowed by having detailed pressure bounds using the phase envelope; thus, smaller diameters can be used. In addition to the heat exchange with the surrounding soil, we also incorporate temperature changes caused by pressure drops in pipelines, namely the~ \textit{Joule-Thomson~effect}. To the best of our knowledge, we are the first to investigate the arc capacity optimization problem in terms of \coTwo transport. 

The remainder of this paper is structured as follows: In Section~\ref{sec:modelling}, we introduce the general problem of discrete pipeline sizing in tree-shaped networks in \ref{netDesign}, followed by the description of the thermophysical effects in the network and pipelines in \ref{sec:thermophysical_modelling}. We present the proposed algorithm that combines these two in Section~\ref{sec:algorithm}. After the model is validated with respect to the pressure and temperature behavior on a single pipeline against data from the literature in Section~\ref{sec:validation}, we apply our approach to a real-world network planning problem in Section~\ref{sec:real_co2_network}. In Section~\ref{sec:robust_analysis}, we introduce an artificially created dataset to test the robustness of our approach.
We conclude the paper with some closing remarks and an outlook in Section~\ref{sec:conclusion}.

\section{Modelling} \label{sec:modelling}
In this section, we first introduce the general network design problem for tree-shaped networks with multiple sources and a single sink. And secondly, we describe the thermophysical effects of transporting fluids in pipeline networks. 

\subsection{Network design problem}\label{netDesign}

We model the \coTwo network as a directed tree-shaped graph $\generalGraph$, where $\setVertices$ and $\setArcs$ represent the set of nodes and the set of arcs in the network, respectively. The set of nodes consists of a set of entry nodes $\setEntries$, a set of exit nodes \setExits, and a set of inner nodes $\setInnerNodes$, i.e., $\setVerticesAll$. Note that we only consider a single exit in this study, i.e.,  $|\setExits| = 1$. We denote the union of the entry nodes and the exit as the set of boundary nodes as $\setBoundaryNodesAll$. The set of arcs $\setArcs$ is distinguished by a set of pipes $\setPipes$ and a set of pumps $\setPumps$. For each pipe $\arc \in \setPipes$, we consider a discrete finite set of possible diameters \setDiameter{\arc}. Moreover, \arcCost{\arc}{\diameter} $\in \Rpos$ describes the related costs of implementing an arc \arc of diameter \diameter. For each pipe $\arc \in \setPipes$ and each diameter $\diameter \in \setDiameter{\arc}$, we introduce a binary decision variable $x_{a,d}$, whether diameter \diameter is chosen for pipe \arc ($x_{a,d}=1$) or not ($x_{a,d}=0$). We denote the set of incoming and outgoing arcs of a node $\vertex \in \setVertices$ as 
\begin{align*}
    &\setIncomingArcsDef{\vertex}{\anotherVertex} &\text{and} &&\setOutgoingArcsDef{\vertex}{\anotherVertex},
\end{align*}
respectively. For each $\vertex \in \setEntries$, we have $\setIncomingArcs{\vertex} = \emptyset$ and $|\setOutgoingArcs{\vertex}| = 1$ and for the exit node $\vertex \in \setEntries$ $|\setIncomingArcs{\vertex}| > 0$ and $\setOutgoingArcs{\vertex} = \emptyset$, i.e., the entry nodes are the leaves of an in-tree rooted at the exit node. Further, let $\setMiddleNodes \subseteq \setInnerNodes$ be the set of intermediate nodes for which we have exactly one incoming and one outgoing pipe with at least two choosable diameters:
\begin{align*}
    \setMiddleNodes \coloneqq \lbrace \vertex \in \setVertices \,|\,  |\setIncomingArcs{\vertex}| = |\setOutgoingArcs{\vertex}| = 1 \land \setIncomingArcs{\vertex}\cup \setOutgoingArcs{\vertex} \in \setPipes \land |\setDiameter{\arc=(\anotherVertex, \vertex)}| > 1 \land |\setDiameter{\arc=(\vertex, w)}| > 1 \rbrace.
\end{align*}
For each node $\vertex\in\setVertices$, a boundary value $\boundaryVal{\vertex}\in \mathbb{R}$ is given representing its inflow ($\boundaryVal{\vertex}<0$) or its outflow ($\boundaryVal{\vertex}<0$). For each $\vertex\in\setEntries$, we have $\boundaryVal{\vertex} \in \mathbb{R}_{\ge 0}$, for each $\vertex\in\setExits$, we have $\boundaryVal{\vertex} \in \mathbb{R}_{\le 0}$, and for each $\vertex \in \setInnerNodes$, we assume $\boundaryVal{\vertex}=0$. Additionally, we have a balanced network, i.e., $\sum_{\vertex\in\setBoundaryNodes} \boundaryVal{\vertex} = 0$. The flow over an arc $\arc\in\setArcs$ is denoted as $\arcFlow{\arc} \in \mathbb{R}_{\ge 0}$. At each node, the flow balance equation has to hold, i.e.,

\begin{align}
& \sum_{\arc\in\setOutgoingArcs{\vertex}}\arcFlow{\arc} - \sum_{\arc\in\setIncomingArcs{\vertex}}\arcFlow{\arc} = \boundaryVal{\vertex} & \forall \vertex \in \setVertices.
\end{align}

For each node $\vertex \in \setVertices$, we introduce a variable representing its pressure level $\vertexPressure{\vertex}\in [\LB{\pressure}{\vertex}, \UB{\pressure}{\vertex}] \subseteq\Rpos$, where $\LB{\pressure}{\vertex}, \UB{\pressure}{\vertex}$ represent the lower and the upper bound, respectively. In general, flow goes from higher pressure levels to lower pressure levels. Pumps are used to increase the pressure of the fluid. For a pipe $\arc=(\vertex, \anotherVertex)\in\setPipes$, the pressure difference between node \vertex and node \anotherVertex is described by the flow \arcFlow{\arc}, the geographical height of the nodes, and the friction loss coefficient $\pipeLossCoeff_\arc(\pipeDiameter{\arc}, \arcFlow{\arc})$ of the pipeline, which depends on the variable pipe diameter \pipeDiameter{\arc} and the flow \arcFlow{\arc}. 
The pressure difference is modeled as
\begin{align}
&\vertexPressure{\vertex} - \vertexPressure{\anotherVertex} = \pipeLossCoeff_\arc\!\left(\pipeDiameter{\arc},\arcFlow{\arc}\right) \arcFlow{\arc} |\arcFlow{\arc}| - (H_\vertex^0 - H_\anotherVertex^0) \rho g& \forall \arc=(\vertex, \anotherVertex)\in\setPipes
\end{align}
where $H_\vertex^0$, $H_\anotherVertex^0$, $\rho$, $g$ denote the geographical height of node \vertex and \anotherVertex, the density of the fluid, and the gravitational constant, respectively. 
The friction loss coefficient $\pipeLossCoeff_\arc$ of a pipe $a$ is modelled using the Weymouth equation
\begin{align}
    \pipeLossCoeff_\arc\!\left(\pipeDiameter{\arc},\arcFlow{\arc}\right) = \frac{8\arcLength}{\pi^2\density\pipeDiameter{\arc}^5}\frictionFactor_\arc(\pipeDiameter{\arc}, \arcFlow{\arc})\label{eq:loss_coeff}
\end{align}
where \arcLength represents the length of the pipe and $\frictionFactor_\arc$ denotes the Darcy-Weisbach friction factor, see \cite{lurie2008}. In this study, we use the implicit Colebrook-White equation, which is the most accurate approximation in case of turbulent flow, see \cite{fugenschuh2015}. The Colebrook-White equation is given by
\begin{align}
    \frac{1}{\sqrt{\frictionFactor_\arc}} = -2 \log_{10} \left(\frac{\varepsilon}{3.7 \pipeDiameter{\arc}} + \frac{2.51}{\Rey\sqrt{\frictionFactor_\arc}}\right) \label{colebrook}
\end{align}
where the roughness of the pipe is denote by $\varepsilon$ and the Reynolds number is denoted by \Rey. The dimensionless Reynolds number is defined as
\begin{align}
&\Rey = \frac{\diameter_\arc\arcFlow{\arc}}{A_\arc \dynViscosity}\label{reynolds}
\end{align}
where $A_\arc$ and $\dynViscosity$ represent the cross area of the pipe and the dynamic viscosity of the fluid, respectively.

The pressure increase generated by a pump $\arc=(\vertex, \anotherVertex) \in \setPumps$ is modelled as follows
\begin{align}
    &\vertexPressure{\anotherVertex}\ge \vertexPressure{\vertex} &\forall \arc=(\vertex, \anotherVertex) \in \setPumps \\
    &\LB{\pressure}{\arc}\le\vertexPressure{\vertex}, \vertexPressure{\anotherVertex}\le \UB{\pressure}{\arc}&\forall \arc=(\vertex, \anotherVertex) \in \setPumps
\end{align}
where \LB{\pressure}{\arc} and \UB{\pressure}{\arc} denote the technical pressure range of the pump $\arc$.

The goal of this network design problem is to find the pipe diameter \pipeDiameter{\arc} for each $\arc \in \setPipes$ of minimal costs such that there exists a feasible state of the network, i.e., all node balances hold and the pressure limits are not violated. The problem $\mathcal{P}$ reads as follows
\begin{align}
    \min_{x,\massFlow, \pressure}& \sum_{\arc\in\setPipes}\sum_{\diameter\in\setDiameter{\arc}} \arcCost{\arc}{\diameter} x_{\arc,\diameter}&\label{diameterOpt}\\
    \text{s.t.}&&\nonumber\\
    & \sum_{\vertex\in\setBoundaryNodes} \boundaryVal{\vertex} = 0,& \label{netBalance}\\
    & \sum_{\arc\in\setOutgoingArcs{\vertex}}\arcFlow{\arc} - \sum_{\arc\in\setIncomingArcs{\vertex}}\arcFlow{\arc} = \boundaryVal{\vertex} & \forall \vertex \in \setVertices,\label{nodebalance}\\
    &\vertexPressure{\vertex} - \vertexPressure{\anotherVertex} + (H_\vertex^0 - H_\anotherVertex^0) \rho g = \sum_{\diameter\in\setDiameter{\arc}} x_{\arc,\diameter}\pipeLossCoeff_\arc\!\left(\diameter,\arcFlow{\arc}\right) \arcFlow{\arc} |\arcFlow{\arc}|& \forall \arc=(\vertex, \anotherVertex)\in\setPipes,\label{pressurePipe}\\
    &\vertexPressure{\anotherVertex}\ge \vertexPressure{\vertex} &\forall \arc=(\vertex, \anotherVertex) \in \setPumps \label{pumpIncrease}\\
    &\LB{\pressure}{\arc}\le\vertexPressure{\vertex}, \vertexPressure{\anotherVertex}\le \UB{\pressure}{\arc}&\forall \arc=(\vertex, \anotherVertex) \in \setPumps,\label{pumpRange}\\
    &\LB{\pressure}{\vertex}\le \vertexPressure{\vertex} \le \UB{\pressure}{\vertex}&\forall \vertex\in\setVertices,\label{pressureLimit}\\
    &\sum_{\diameter\in\setDiameter{\arc}}x_{\arc,\diameter} = 1 & \forall \arc \in \setPipes,\label{arcSizeDecision}\\
    &\sum_{\diameter\in\setDiameter{(\anotherVertex, \vertex)}}x_{(\anotherVertex, \vertex),\diameter}\,\diameter = \sum_{\diameter\in\setDiameter{(\vertex, w)}}x_{(\vertex, w),\diameter}\,\diameter & \forall \vertex \in \setMiddleNodes,\label{sameDiameter}\\
    & x_{\arc,\diameter}\in \lbrace0,1\rbrace&\forall \arc \in \setPipes, \forall \diameter \in \setDiameter{\arc}.\label{x_diameter}
\end{align}
This model incorporates the balance of the network \eqref{netBalance}, the node balance \eqref{nodebalance}, the pressure drop along a pipeline for a chosen diameter \diameter \eqref{pressurePipe}, the pressure increase and range of pumps \eqref{pumpIncrease} and \eqref{pumpRange}, the pressure limits on the nodes \eqref{pressureLimit}, the choice of exactly one diameter for each pipeline \eqref{arcSizeDecision}, and that for the incident pipelines of node $\vertex \in \setMiddleNodes$ the same diameter is chosen \eqref{sameDiameter}. Note that there may exist pipelines in the network, which have a prefixed diameter, and, thus, \eqref{sameDiameter} does not apply and a different diameter is allowed on the previous and the next pipeline.

In general, $\mathcal{P}$ is a mixed integer (discrete decisions regarding pipe diameter) nonlinear (e.g., the behavior of the friction loss coefficient) mathematical program (MINLP). Additionally, as mentioned in Section~\ref{sec:introduction}, physical properties are nonlinear functions of pressure \pressure and temperature \temperature, e.g., the density \density is defined as $\density = \density(\pressure,\temperature)$. However, the complexity of the problem can be drastically reduced by exploiting the fact that we have a tree-shaped graph rooted at a single exit node with all the entry nodes at the leaves and all arcs directed towards the root. Thus, for given boundary values $\boundaryVal{\vertex}$ for all $\vertex \in \setEntries$, we compute the unique values of the flow variables and friction loss coefficient with respect to flow and each diameter a priori, which is described in Section~\ref{sec:algorithm}. 

\subsection{Thermophysical modelling}\label{sec:thermophysical_modelling}
In the following, we describe three main effects that impact the temperature changes of the flow in the network: Namely, the temperature of fluid mixtures at junction points in the network, heat exchange with the surrounding alongside the pipeline, and the \textit{Joule-Thomson effect} which describes the temperature change with respect to pressure change. 
\subsubsection{Temperature of a mixture}
The fluid temperature $\temperature_\vertex$ at a junction point $\vertex$ is derived from the conservation of energy and is determined by
\begin{align}
    \temperature_\vertex = \frac{\sum_{\arc \in \setIncomingArcs{\vertex}}\heatCapArc{\arc}\arcFlow{a}\outTempArc{\arc}}{\sum_{\arc \in \setIncomingArcs{\vertex}}\heatCapArc{\arc}\arcFlow{a}} \label{mixingTemp}
\end{align}
where $\heatCapArc{\arc}$, $\arcFlow{\arc}$, and $\outTempArc{\arc}$ represent the heat capacity, the mass flow, and the outlet temperature of the incoming arc $\arc$, respectively, see e.g. \cite{fugenschuh2015, schmidt2015}.
For the remainder of this section the index $a$ on pipeline specific quantities is dropped for the sake of clarity. 
\subsubsection{Energy balance of a pipeline}
Under the assumption of a steady-state, one-dimensional flow in the pipeline, by neglecting the change in the fluid kinetic energy and conduction, the energy balance across a fluid segment of length $dx$ is
\begin{align}
    \massFlow\frac{d\enthalpy}{dx} = - \heatFlow_w \label{energyBalance}
\end{align}
where $\heatFlow_w$ is the heat loss per unit length, which consists of the heat exchange with the the surrounding soil $\heatFlow_c$, and the friction heat loss occurring at the wall $\heatFlow_f$, i.e.,
\begin{align}
    \heatFlow_w = \heatFlow_c + \heatFlow_f. \label{heat_loss}
\end{align}
Under the assumption that \enthalpy is a function of pressure and temperature, i.e., $\enthalpy = \enthalpy(\pressure, \temperature)$, we can express $d\enthalpy$ as
\begin{align}
    d\enthalpy = \heatCap d\temperature - \heatCap\jouleThomson d\pressure, \label{jouleThomson}
\end{align}
where \jouleThomson is the \textit{Joule-Thomson} coefficient \cite{colina1998}. 
Applying \eqref{heat_loss} and \eqref{jouleThomson} on \eqref{energyBalance}, dividing by \massFlow and \heatCap, and assuming $\heatFlow_f\ll \heatFlow_c$ and, thus, neglecting $\heatFlow_f$, gives
\begin{align}
    \frac{d\temperature}{dx} = \jouleThomson \frac{d\pressure}{dx} - \frac{\heatFlow_c}{\massFlow\heatCap} \label{tempChange_dx} 
\end{align}
where the first term on the right hand side represents the \textit{Joule-Thomson} effect due to pressure changes and the second term the heat exchange with the surrounding. Integrating \eqref{tempChange_dx} over the length of the pipeline \length gives
\begin{align}
    \outTemp - \inTemp = \jouleThomson (\outPressure - \inPressure) - \frac{\length}{\massFlow\heatCap}\heatFlow_c. \label{tempChange}
\end{align}
For a detailed description the reader is referred to \cite{phuoc2019}. Next, the heat exchange with the surrounding is described in more detail.
\subsubsection{Heat transfer with surrounding}
The heat exchange $\Delta \heatFlow$ with the surrounding along a pipeline $a\in \setPipes$ is calculated using the correlation for turbulent throughflow of a buried pipeline. A detailed description can be found in Chapter E and G of \cite{VDIWaerme2013}. For the heat exchange $\Delta \heatFlow$ we have
\begin{align}
    \Delta \heatFlow = \length \heatFlow_c = \heatTransRate \length \logTemp \label{heatExchange}
\end{align}
where \heatTransRate denotes the heat transmission rate over the pipeline and \logTemp is the logarithmic temperature difference given by
\begin{align}
    \logTemp = \frac{\inTemp - \outTemp}{\log{\frac{\inTemp - \soilTemp}{\outTemp- \soilTemp}}}.
\end{align}
The heat transmission rate \heatTransRate consists of the heat transmission through the fluid inside the pipeline, the transmission through the pipeline wall, and finally, the heat transmission into the surrounding soil as we consider pipelines buried in the ground. The heat transmission rate is defined as
\begin{align}
    \heatTransRate = \frac{2\pi}{\frac{2}{\heatTransCoeff \diameter} + \frac{1}{\thermCond_{\text{p}}} \log{\frac{\outerDiameter}{\diameter}} + \frac{1}{\thermCond_{\text{s}}} \log{\frac{4\pipeDepth}{\diameter}}}
\end{align}
where $\thermCond_{\text{p}}$, $\thermCond_{\text{s}}$, \pipeDepth, and \outerDiameter represent the thermal conductivity of the pipeline material and the surrounding soil, the buried depth, and the outer diameter of the pipeline, respectively. Whereas the heat transmission through the pipeline wall into the surrounding soil is mainly dependent on given material parameters, the heat transmission factor of the fluid, denoted by \heatTransCoeff, also depends on the flow. To determine \heatTransCoeff, we introduce the following well-known dimensionless quantities: The \textit{Prandtl} number, which is defined as the ratio of the dynamic viscosity to the thermal diffusivity as
\begin{align}
    \Pra = \Pra(\dynViscosity, \thermCond_{\text{f}}, \heatCap) = \frac{\dynViscosity}{\thermCond_{\text{f}}/\heatCap} = \frac{\heatCap \dynViscosity}{\thermCond_{\text{f}}}, \label{prandtl}
\end{align}
where $\thermCond_{\text{f}}$ is the thermal conductivity of the fluid, and the \textit{Nusselt} number, which describes the ratio of the convective and conductive heat transfer as
\begin{align}
    &\Nus = \frac{\heatTransCoeff\diameter}{\thermCond_{\text{f}}} & \Longleftrightarrow && \heatTransCoeff = \frac{\Nus \thermCond_{\text{f}}}{\diameter}.\label{nusselt_a}
\end{align}
For fully developed turbulent flow, the relationship between the heat transfer and flow resistance is given by
\begin{align}
    \Nus = \frac{(\zeta/8)\,\Rey\,\Pra}{1 + 12.7 \sqrt{\zeta/8}(\Pra^{2/3} -1 )} \left[ 1 + \left(\frac{\diameter}{\length}\right)^{2/3}\right] \label{nusselt_b}
\end{align}
with 
\begin{align}
    \zeta = \left(1.8 \log_{10}\Rey -1.5\right)^{-2},
\end{align}
where \Rey is the Reynolds number given in \eqref{reynolds}. Again, for a more detailed description and derivation of these characteristic numbers, we refer the reader to~\cite{VDIWaerme2013}. Using \eqref{reynolds}, \eqref{prandtl}, \eqref{nusselt_a}, and \eqref{nusselt_b} we calculate the heat transmission factor of the fluid \heatTransCoeff.

\subsubsection{Outlet temperature of a pipeline}
By applying \eqref{heatExchange} on \eqref{tempChange} the outlet temperature of a pipeline is given by
\begin{align}
    \outTemp = \inTemp + \jouleThomson(\outPressure - \inPressure) - \frac{1}{\massFlow\heatCap}\heatTransRate \length \logTemp. \label{outTemp}
\end{align}

\section{Algorithm}\label{sec:algorithm}

As mentioned in the introduction, the behavior of \coTwo is very sensitive to temperature and pressure changes, which can significantly change physical properties such as density. 
In order to obtain an optimal solution to the network design problem taking changes in the physical properties of \coTwo due to pressure and temperature changes into account, we propose an iterative algorithm splitting the problem into the network design problem described in \ref{netDesign} and the thermophysical modeling of a network with given pipeline diameters and pressure levels introduced in \ref{sec:thermophysical_modelling}. 
The algorithm consists of four main parts: First, the computation of the unique flow values in the network. Second, the determination of the friction factors for each possible pipeline diameter. Third, computation of the temperature distribution in the network. And finally, solving the network design problem. In the following, these four parts are described in more detail before we formulate the complete algorithm.

\subsection{Computing unique flow values}\label{sec:computing_flow}
Computing the unique flow variables is done as follows. Considering the directed in-tree graph rooted at a single exit node introduced above, let $\setVertices(\anotherVertex) \coloneqq \lbrace \vertex \in \setVertices : (\anotherVertex, \vertex) \in \setArcs \lor (\vertex, \anotherVertex) \in \setArcs\rbrace $ be the set of nodes adjacent to node \anotherVertex. Further, let $L\coloneqq\lbrace \vertex \in \setVertices: |\setVertices(\vertex)|=1 \land \vertex\notin \setExits\rbrace$ be the set of leaves in the tree $\graph$. For each $\vertex \in L$ set the flow of the unique outgoing arc $\arc \in \setOutgoingArcs{\vertex}$ to $\arcFlow{\arc} = \boundaryVal{\vertex}$ and update $\boundaryVal{\anotherVertex} \coloneqq \boundaryVal{\anotherVertex} + \arcFlow{\arc}$ for each neighbor $\anotherVertex \in \setVertices(\vertex)$. Let $\setVertices = \setVertices\backslash\setLeafs$ and $\setArcs = \setArcs \backslash \bigcup_{\vertex\in\setLeafs}\setOutgoingArcs{\vertex}$, update \setLeafs, and then iterate until \setArcs is empty.

\subsection{Computing friction loss coefficient}\label{sec:computing_friction}

For each pipe $\arc\in\setPipes$ and diameter $\diameter \in \setDiameter{\arc}$, the friction loss coefficient is computed using the computed flow values from Section~\ref{sec:computing_flow} and \eqref{eq:loss_coeff}. To determine the friction factor $\frictionFactor_\arc$ in \eqref{eq:loss_coeff}, we solve \eqref{colebrook} iteratively until its change is below a given tolerance \frictionTol.

\subsection{Temperature distribution in the network}\label{sec:computing_temperature}
Similar to the algorithm used for precomputing the unique flow values in a tree-shaped network described in Section~\ref{sec:computing_flow}, the temperature distribution is calculated. Let \generalGraph be a directed in-tree graph as described in \ref{netDesign} with given flow value $\arcFlow{\arc}$ for each $\arc \in \setArcs$ and pressure level $\vertexPressure{\vertex}$ for each $\vertex\in \setVertices$. Again, let $L\coloneqq\lbrace \vertex \in \setVertices: |\setVertices(\vertex)|=1 \land \vertex\notin \setExits\rbrace$ be the set of leafs in the tree $\graph$. For each $\vertex \in \setEntries$, an inflow temperature $\temperature_\vertex = \inTemp$ is given. The temperature distribution is computed as follows: For each $\vertex \in \setLeafs$, set the inlet temperature of the unique outgoing arc $\arc \in \setOutgoingArcs{\vertex}$ to $\inTempArc{\arc} = \temperature_\vertex$. If $\arc \in \setPipes$, the outlet temperature \outTempArc{\arc} is determined by solving \eqref{outTemp} iteratively. For a given pressure difference $\Delta\pressure = \outPressure - \inPressure$ and inlet temperature of the pipeline, starting with an initial guess $\outTemp{}^{,0}$, the outlet temperature of the $i$-th iteration is calculated using physical properties at the mean temperature $\temperature_\mathrm{m} = (\inTemp + \outTemp{}^{,i-1})/2$ and pressure $\pressure_\mathrm{m} = (\inPressure + \outPressure) /2$, i.e.,
\begin{align*}
    \heatCap(\pressure_m,\temperature_m),\,\dynViscosity(\pressure_m,\temperature_m),\,\thermCond_{\text{f}}(\pressure_m,\temperature_m),\,\jouleThomson(\pressure_m,\temperature_m).
\end{align*}
We iterate until $|\outTemp{}^{,i} -\outTemp{}^{,i-1}| < \temperatureTol$. If $\arc \in \setPumps$, set $\outTempArc{\arc} = \inTempArc{\arc}$. Then, for each node $\anotherVertex \in \bigcup_{\vertex\in\setLeafs}\setVertices(\vertex)$, compute its mixing temperature $\temperature_\anotherVertex$ using \eqref{mixingTemp}, where $\setVertices(\vertex)$ denotes the set of neighbors of \vertex. Again, let $\setVertices = \setVertices\backslash\setLeafs$ and $\setArcs = \setArcs \backslash \bigcup_{\vertex\in\setLeafs}\setOutgoingArcs{\vertex}$, update \setLeafs, and iterate until \setArcs is empty. 

\subsection{Solving the network design problem}

Using the flow values computed in Section~\ref{sec:computing_flow} and the friction loss coefficient from Section~\ref{sec:computing_friction}, the network design problem in \eqref{diameterOpt}--\eqref{x_diameter} is solved as a tractable mixed integer linear program.

\subsection{Complete algorithm}\label{sec:complete_algo}

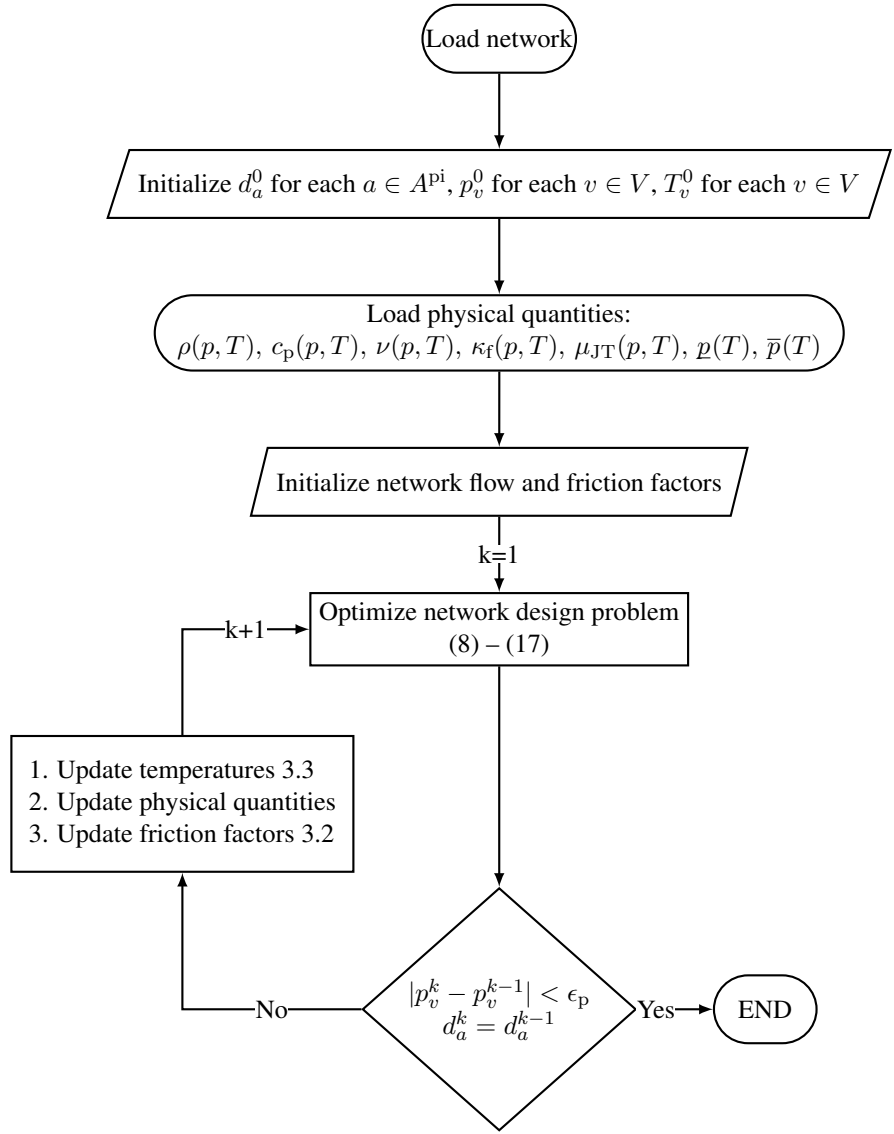
\begin{figure}
    \centering
    \begin{tikzpicture}[scale=0.9, every node/.style={scale=0.9}, font=\normalsize,thick]
	\node[draw,
	rounded rectangle,
	minimum width=2.5cm,
	minimum height=1cm] (block1) {Load network};
	\node[draw,
	trapezium, 
	trapezium left angle = 65,
	trapezium right angle = 115,
	trapezium stretches,
	below=of block1,
	minimum width=3.5cm,
	minimum height=1cm
	] (block2a) { Initialize $\pipeDiameter{\arc}^0$ for each $\arc \in \setPipes$, $\vertexPressure{\vertex}^0$ for each $\vertex\in\setVertices$, $\temperature_\vertex^0$ for each $\vertex\in\setVertices$};
	
	\node[draw,
	rounded rectangle,
	below=of block2a,
	minimum width=2.5cm,
	minimum height=1cm,
	align=center,
	] (block2b) {Load physical quantities: \\  $\density(\pressure,\temperature),\,\heatCap(\pressure,\temperature),\,\dynViscosity(\pressure,\temperature),\,\thermCond_{\text{f}}(\pressure,\temperature),\,\jouleThomson(\pressure,\temperature),\, \LB{\pressure}{}(\temperature),\, \UB{\pressure}{}(\temperature)$};
	
	\node[draw,
	trapezium, 
	trapezium left angle = 65,
	trapezium right angle = 115,
	trapezium stretches,
	below=of block2b,
	minimum width=3.5cm,
	minimum height=1cm
	] (block2) { Initialize network flow and friction factors};
	
	\node[draw,
	below=of block2,
	minimum width=3.5cm,
	minimum height=1cm,
	align=center,
	] (block3) {Optimize network design problem \\ \eqref{diameterOpt} -- \eqref{x_diameter}};
	
	\node[
	below=of block3,
	minimum height=1cm,
	](blockDummy) {};
	
	\node[draw,
	diamond,
	below=of blockDummy,
	minimum width=4.cm,
	inner sep=0,
	align=center] (block4) { $|\vertexPressure{\vertex}^{k} - \vertexPressure{\vertex}^{k-1}|<\pressureTol $\\ $\pipeDiameter{\arc}^{k} =  \pipeDiameter{\arc}^{k-1}$};
	
	\node[draw,
	rounded rectangle,
	right=of block4,
	minimum width=1.5cm,
	minimum height=1cm,
	inner sep=0] (block5) { END };
	
	\node[draw,
	above left=of block4,
	minimum width=5cm,
	minimum height=2cm,
	align=left,
	inner sep=0] (block6) { 1. Update temperatures \ref{sec:computing_temperature} \\2. Update physical quantities\\ 3. Update friction factors \ref{sec:computing_friction}};	
	
	
	\draw[-latex] (block1) edge (block2a);
	\draw[-latex] (block2a) edge (block2b);
	\draw[-latex] (block2b) edge (block2);
	\draw[-latex] (block2) -- (block3)
	node[pos=0.5,fill=white,inner sep=0]{k=1};
	\draw[-latex](block3) edge (block4);

	\draw[-latex] (block4) -- (block5)
	node[pos=0.25,fill=white,inner sep=0]{Yes};
	
	\draw[-latex] (block4) -| (block6)
	node[pos=0.25,fill=white,inner sep=0]{No};
	
	
	\draw[-latex] (block6) |- (block3)
	node[pos=0.75,fill=white,inner sep=0]{k+1};	
	\end{tikzpicture}
    \caption{Algorithm to fully solve the thermo-dependent network design problem}
    \label{fig:algo}
\end{figure}

The complete algorithm is structured as follows. In the initial phase, we first choose an initial diameter $\pipeDiameter{\arc}^0$ for each pipeline $\arc \in \setPipes$ and initial pressure level $\vertexPressure{\vertex}^0$ and temperature value $\temperature_\vertex^0$ for each vertex $\vertex\in \setVertices$. Then, we compute the unique flow values for each $\arc\in\setArcs$ and the friction loss coefficient $\pipeLossCoeff_{\arc}(\pipeDiameter{\arc}, \arcFlow{\arc})$ for each pipeline $\arc \in \setPipes$ and possible diameter $\pipeDiameter{\arc}\in \setDiameter{\arc}$. In each iteration, physical properties of \coTwo, such as density, dynamic viscosity, et cetera, are drawn from look-up tables with respect to temperature and pressure levels from the previous iteration. For each $\vertex\in\setVertices$ limits on pressure levels \LB{\pressure}{\vertex} and \UB{\pressure}{\vertex} are updated with respect to the saturation pressure dependent on its temperature. Then, the network design problem $\mathcal{P}$ is solved. The new temperature distribution is computed with the resulting diameters and pressure levels. With the current pressure and temperature levels, the physical properties are updated from the look-up tables and the new friction factors are computed. The algorithm stops if the change in pressure is below the given tolerance \pressureTol and the diameters stay the same.
The full algorithm is shown in Fig. \ref{fig:algo}.

\section{Computational results}\label{sec:computational_results}

In this section, we first test our implementation with respect to the behavior of the temperature and pressure on a single pipeline with a fixed diameter. The results are compared to the pipeline investigated in \cite{nimtz2010}. Then, we apply our proposed algorithm from \ref{sec:algorithm} to a real-world \coTwo network planning problem. Finally, we test the robustness of the algorithm by solving a set of artificially created network instances. The algorithm has been implemented using \python~\oldstylenums{3.8.10}, the network design problem has been implemented using the \gurobi~\python-interface and is solved with \gurobiVersion{9}{5}~\cite{gurobi95}. The computations in Section~\ref{sec:validation} and \ref{sec:real_co2_network} have been performed on an Intel(R) Core(TM) i7-9700K CPU @ 3.60GHz with 64~GB of RAM. We use as tolerances 
\begin{align*}
    &\frictionTol = 10^{-6}, & \temperatureTol = \SI{0.1}{\kelvin}, && \pressureTol = \SI{0.1}{\bar}.
\end{align*}

\subsection{Validation of a single pipeline}\label{sec:validation}
We compare our implementation on a single pipeline taken from \cite{nimtz2010}. The data is given in Table \ref{tab:nimtz_pipe} and Table \ref{tab:boundaries}. Missing information regarding parameters, such as, thermal conductivity of the pipeline and the soil, are assumed. 

\begin{table}[htb]
\centering
\begin{threeparttable}
    \caption{Pipeline data}\label{tab:nimtz_pipe}
    \centering
    \begin{tabular}{lcc}
         Property& Value &Unit\\
         \hline
         Pipe length        & 150   & \si{\kilo\meter}\\
         Pipe diameter      & 0.5   &\si{\meter}\\
         Pipe roughness     & 0.0005&\si{\meter}\\
         Pipe thermal conductivity\tnote{a} & 30.0 &\si[per-mode = symbol]{\watt\per\meter\per\kelvin}
    \end{tabular}
    \begin{tablenotes}
       \item [a] \small Steel pipeline~\cite{nimtz2016}.
     \end{tablenotes}
    \end{threeparttable}
\end{table}

\begin{table}[htb]
\centering
\begin{threeparttable}
    \caption{Boundary conditions}
    \label{tab:boundaries}
    \centering
    \begin{tabular}{p{0.35\textwidth}p{0.1\textwidth}p{0.1\textwidth}}
         Boundary Condition& Value &Unit\\
         \hline
         Mass flow $\arcFlow{\arc}$           & 117   &\si[per-mode = symbol]{\kilogram\per\second}\\
         Inlet temperature          & 40    &\si{\celsius}\\
         Critical Temperature of \coTwo $\temperature_c$ & 31 &\si{\celsius}\\
         Critical pressure of \coTwo $\pressure_c$ & 73.773 &\si{\Bar}\\
         Maximum inlet pressure \UB{\pressure}{\vertexNumbered{\text{in}}}& 150&\si{\Bar}\\
         Minimum outlet pressure \LB{\pressure}{\vertexNumbered{\text{out}}}& 85&\si{\Bar}\\
         Soil temperature           & 10   &\si{\celsius}\\
         Soil thermal conductivity\tnote{b}  & 1.0 &\si[per-mode = symbol]{\watt\per\meter\per\kelvin}
    \end{tabular}
    \begin{tablenotes}
       \item [b] \small Average value from \cite{nimtz2016}.
     \end{tablenotes}
    \end{threeparttable}
\end{table}

To test our implementation on a single horizontal pipeline with fixed diameter $\diameter_{\text{fix}}$, we need to change the network design problem. Since pressure and temperature levels are only computed at the endpoints of a pipeline, we model the pipeline as a tree network with multiple pipe segments \setPipes, a single entry node \vertexNumbered{\text{in}}, and single exit \vertexNumbered{\text{out}} to gain a higher resolution of the temperature and pressure profiles in the pipeline. We split the \SI{150}{\kilo\meter} pipeline into 300 segments of \SI{500}{\meter} each, resulting in $|\setInnerNodes| = |\setMiddleNodes| = 299$ inner nodes. Finally, instead of minimizing the pipeline costs, we minimize the inlet pressure $\vertexPressure{\vertexNumbered{\text{in}}}$. The adapted network design problem from \ref{netDesign} reads as follows

\begin{align}
    \min_{\pressure}& \quad \vertexPressure{\vertexNumbered{\text{in}}}\\
    \text{s.t.}&&\nonumber\\
    &\vertexPressure{\vertex} - \vertexPressure{\anotherVertex} = \pipeLossCoeff_\arc\!\left(\diameter_{\text{fix}},\arcFlow{\arc}\right) \arcFlow{\arc} |\arcFlow{\arc}|& \forall \arc=(\vertex, \anotherVertex)\in\setPipes,\\
    &\vertexPressure{\vertexNumbered{\text{in}}} \le \UB{\pressure}{\vertexNumbered{\text{in}}} ,\,\LB{\pressure}{\vertexNumbered{\text{out}}} \le \vertexPressure{\vertexNumbered{\text{out}}}&
\end{align}

\begin{figure}[htb]
\centering
\begin{subfigure}[b]{.45\linewidth}
\includegraphics[width=\linewidth]{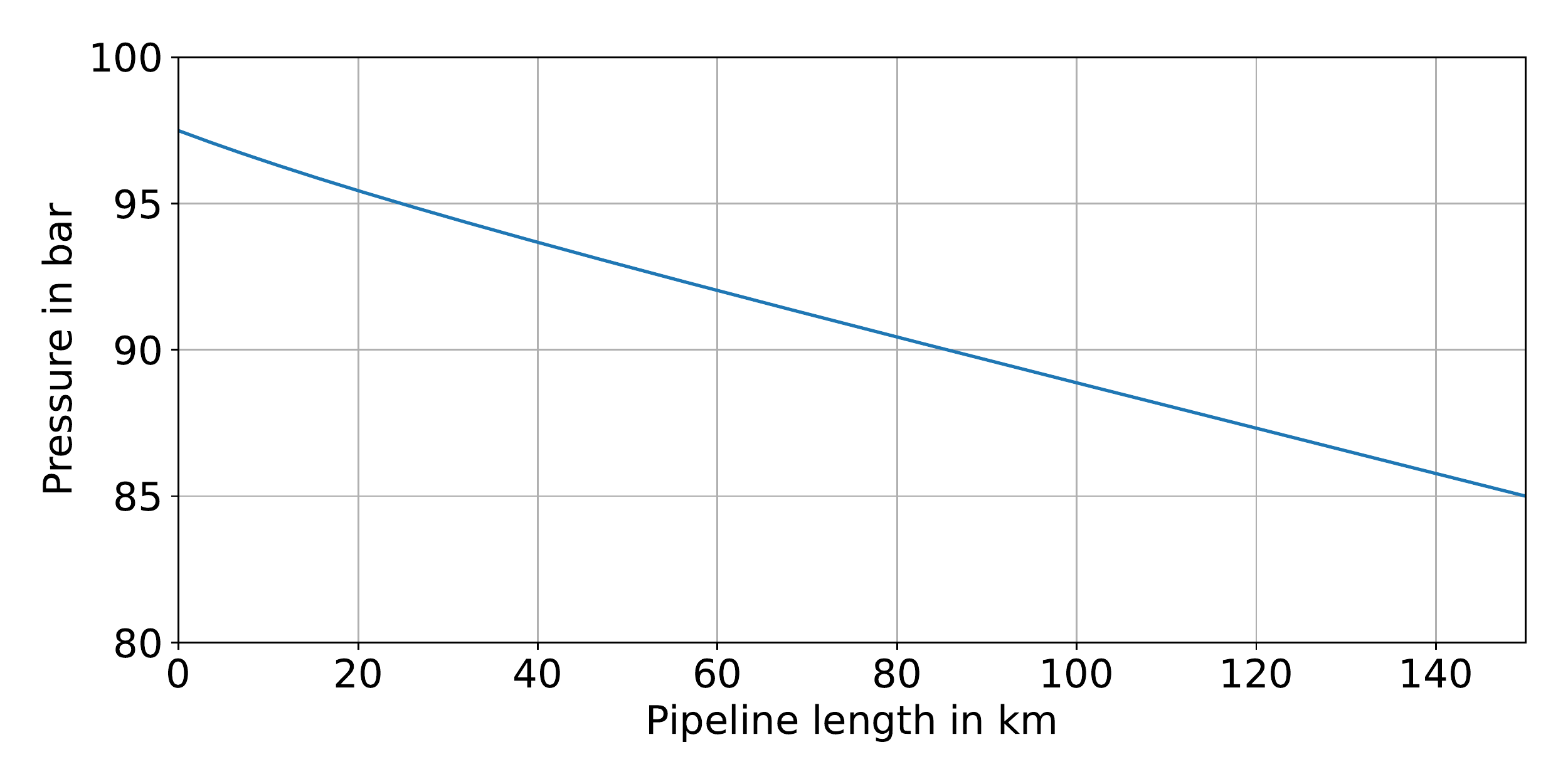}
\caption{Pressure}\label{fig:mouse}
\end{subfigure}
\begin{subfigure}[b]{.45\linewidth}
\includegraphics[width=\linewidth]{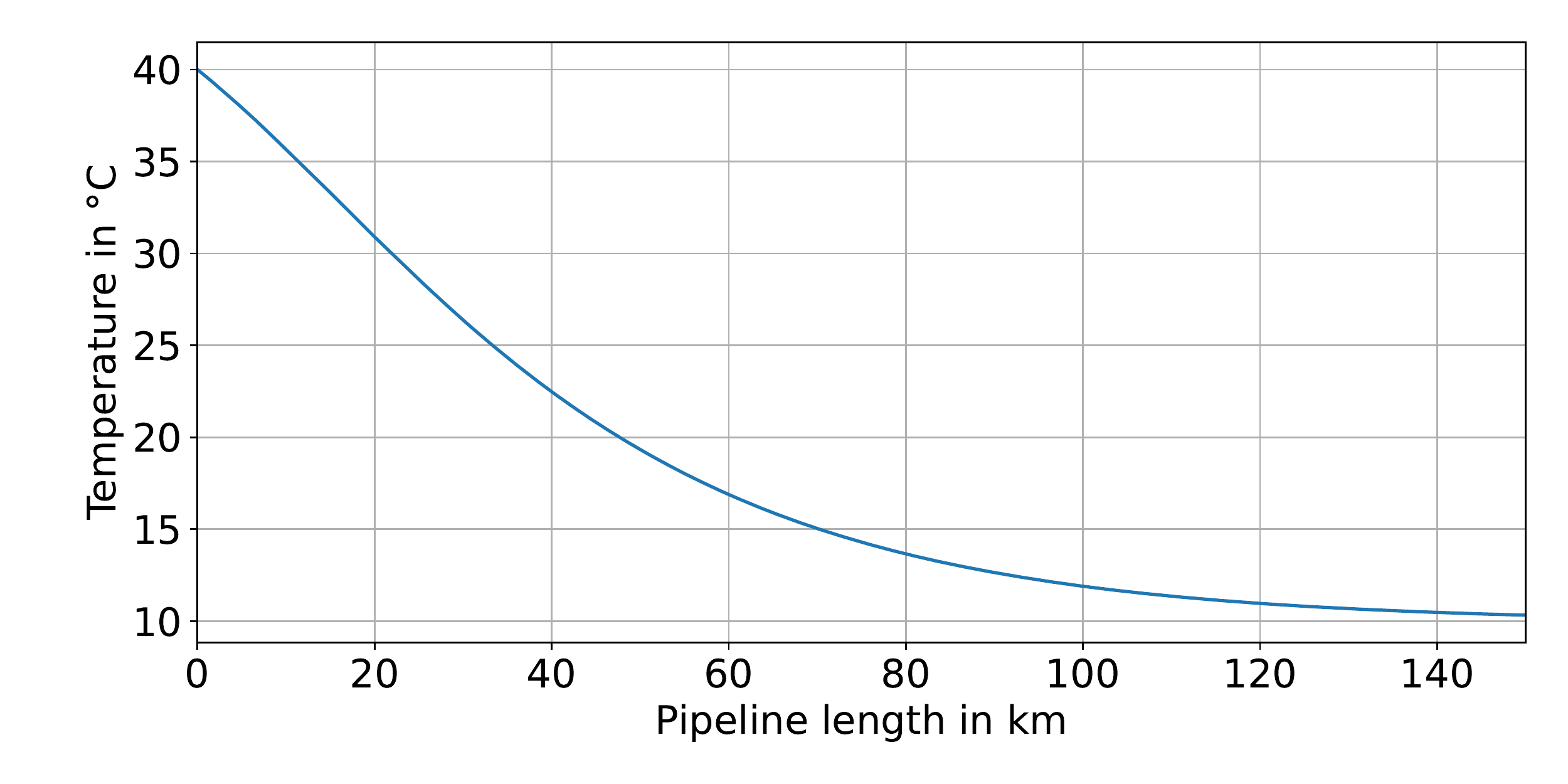}
\caption{Temperature}\label{fig:gull}
\end{subfigure}
\begin{subfigure}[b]{.45\linewidth}
\includegraphics[width=\linewidth]{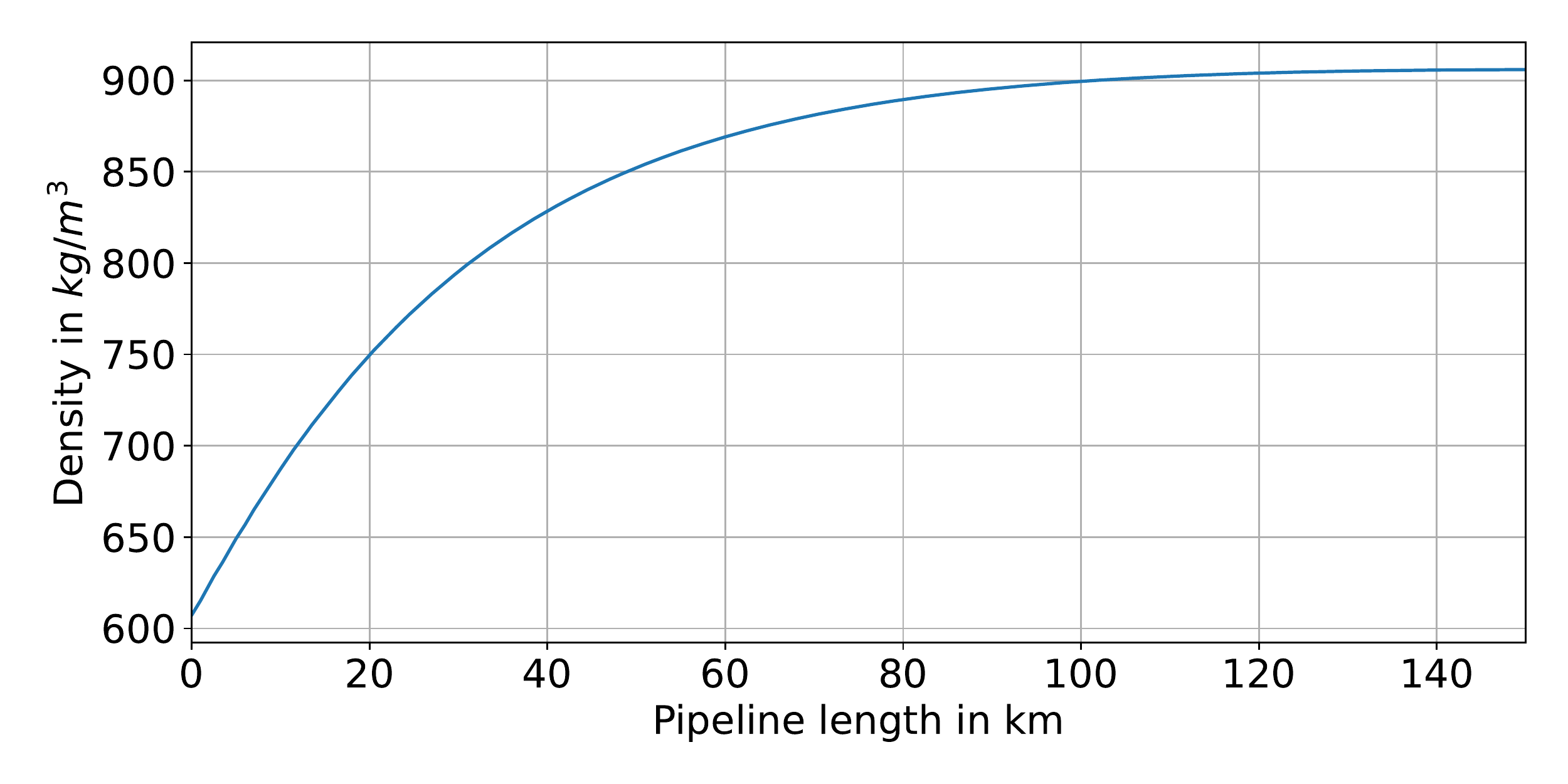}
\caption{Density}\label{fig:tiger}
\end{subfigure}
\caption{Resulting \coTwo behavior along the pipeline for validation}
\label{fig:results_validate_nimtz2010}
\end{figure}

The results are given in Fig. \ref{fig:results_validate_nimtz2010}. The behavior of the \coTwo is inline to the results given for \textit{case~1} in \cite{nimtz2010}. The inlet pressure is at \SI{97.5}{\Bar} and the outlet temperature at \SI{10.33}{\celsius}, i.e., \SI{283.48}{\kelvin}, compared to \SI{98}{\Bar} and \SI{12}{\celsius}, i.e., \SI{285.15}{\kelvin}, resulting in a deviation of less than one percent. The outlet pressure has with \SI{85}{\bar} the same value. A comparison to the other cases given in~\cite{nimtz2010} is not meaningful, since we do not aim for a model including a storage well in our study.

\subsection{Real-world CO2 transport network}\label{sec:real_co2_network}
In the following, the instance of real-world network planning problem for \coTwo transport and the computational results are presented. 
\subsubsection{Network data}\label{sec:net_data}
\begin{figure}
\centering
\includegraphics[width=0.65\linewidth]{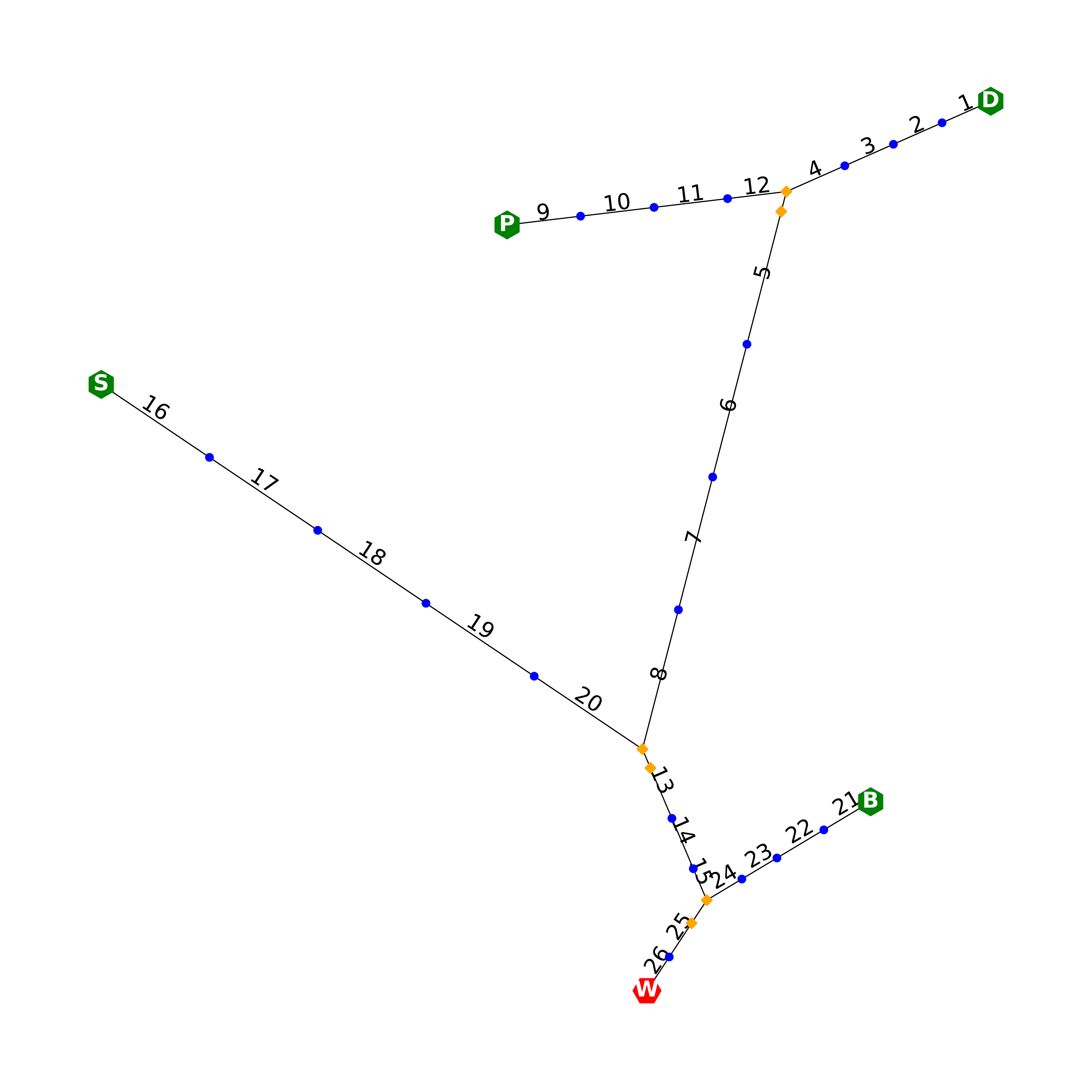}
\caption{Original network topology: Supply nodes S, P, D, B (\markerentry), sink node W (\markerexit), tail and head node of pumps (\markerpump), and inner nodes (\markerinnode).}
\label{fig:network_topo}
\end{figure}

In the initial phase of the network operation, the corresponding instance aims to connect multiple cement production sites in North-Western Germany to the North Sea. The topology, the inflow pattern, and the set of possible pipeline diameters are provided by our project partner Open Grid Europe~\cite{oge}. The topology is shown in Figure~\ref{fig:network_topo}. 

The original network consists of four entry nodes, 25 inner nodes, of which 16 nodes are intermediate nodes, 26 pipelines with a total length of 831~\si{\kilo\meter}, and three pumps. Again, each pipeline is split into segments of 500~\si{\meter}, resulting in 1636 additional inner nodes and 1662 pipe segments in the network. The inflow pattern at the four entry nodes is shown in Table~\ref{tab:network_entries}. 

\begin{table}[tbh]
    \centering
    \begin{tabular}{l|c|c|c|c}
                                & D & P & S & B \\
        \hline
         Inflow in \si[per-mode = symbol]{\kg\per\second} & 16 & 8   & 8 & 8\\
         Temperature in \si{\celsius}                       & 50 & 60  & 20& 30\\
         Max Pressure in \si{\bar}                          & 95 & 105 & 88& 88
    \end{tabular}
    \caption{Inflow pattern of entry nodes}
    \label{tab:network_entries}
\end{table}

The lower pressure limit of the sink node is set to 80~\si{\bar}. In general, the lower bound on the pressure levels is given by the saturation pressure dependent on the temperature with an additional security tolerance of \SI{5}{\bar}, such that the fluid stays in the supercritical or liquid form. This lower bound is updated in each iteration. For each pipeline $\arc~\in~\setPipes$ the diameter is chosen in \si{\milli\meter} out of $\setDiameter{\arc}\!~=~\!\lbrace30, 40, 50, 80, 120, 150, 200, 250, 300, 350, 420, 500\rbrace$. The diameters of pipe 1 and 17 are fixed to 300~\si{\milli\meter}. The surrounding soil temperature is assumed to be 10~\si{\celsius}. For all physical properties of pure \coTwo, lookup tables are generated using the software \gascalcVersion{2.6.1}~\cite{gascalc}.

\subsubsection{Results}\label{sec:net_results}

\begin{figure}[htb]
    \centering
    \begin{tikzpicture}
    \begin{axis} [
        width=\textwidth,
        height=0.35\textwidth,
        ybar, 
        xtick=data,
        ylabel={Diameter in \si{\milli\meter}},
        xlabel={Pipeline},
        ymajorgrids=true,
        grid style=dashed,
        enlargelimits=0.05,]
    \addplot table[x=Pipe, y=Diameter]{\mydata};
    \end{axis}
    \end{tikzpicture}
    \caption{Resulting diameter for each pipeline in the network}
    \label{fig:diameter_sol}
\end{figure}
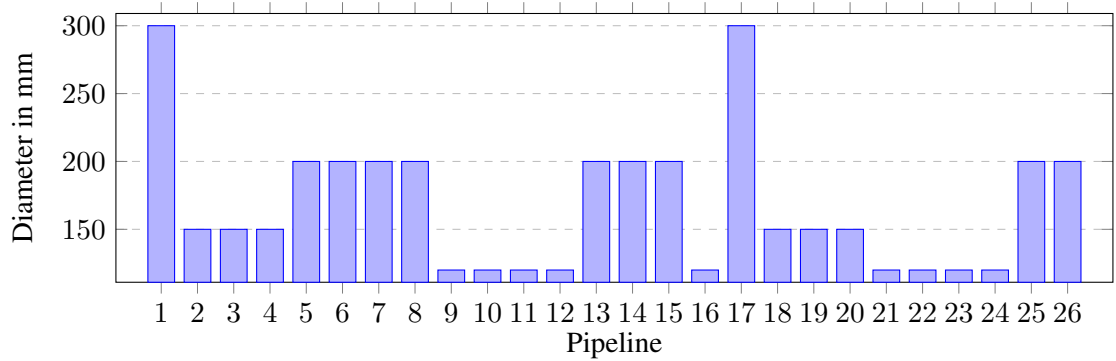

The algorithm stops after five iteration and 28 seconds, the mean time of solving the optimization model is around four seconds, and the time to update the temperature and friction factor is below one second.
The optimal solution consists of \SI{75}{\kilo\meter} pipeline of \SI{300}{\milli\meter} diameter, \SI{292}{\kilo\meter} of \SI{200}{\milli\meter} diameter, \SI{230}{\kilo\meter} of \SI{150}{\milli\meter} diameter, and \SI{234}{\kilo\meter} of \SI{120}{\milli\meter} diameter. The resulting diameters are shown in Figure~\ref{fig:diameter_sol}. 

\pgfplotstableread[row sep=\\, col sep=&]{
Pipe & Diameter\\
1 & 300\\
2 & 150\\
3 & 150\\
4 & 150\\
5& 200\\
6& 200\\
7& 200\\
8& 200\\
9& 120\\
10& 120\\
11& 120\\
12& 120\\
13& 200\\
14& 200\\
15& 200\\
16& 120\\
17& 300\\
18& 150\\
19& 150\\
20& 150\\
21& 120\\
22& 120\\
23& 120\\
24& 120\\
25& 200\\
26& 200\\
}\mydata

\begin{figure}
\centering
\begin{subfigure}[b]{.45\linewidth}
\includegraphics[width=\linewidth]{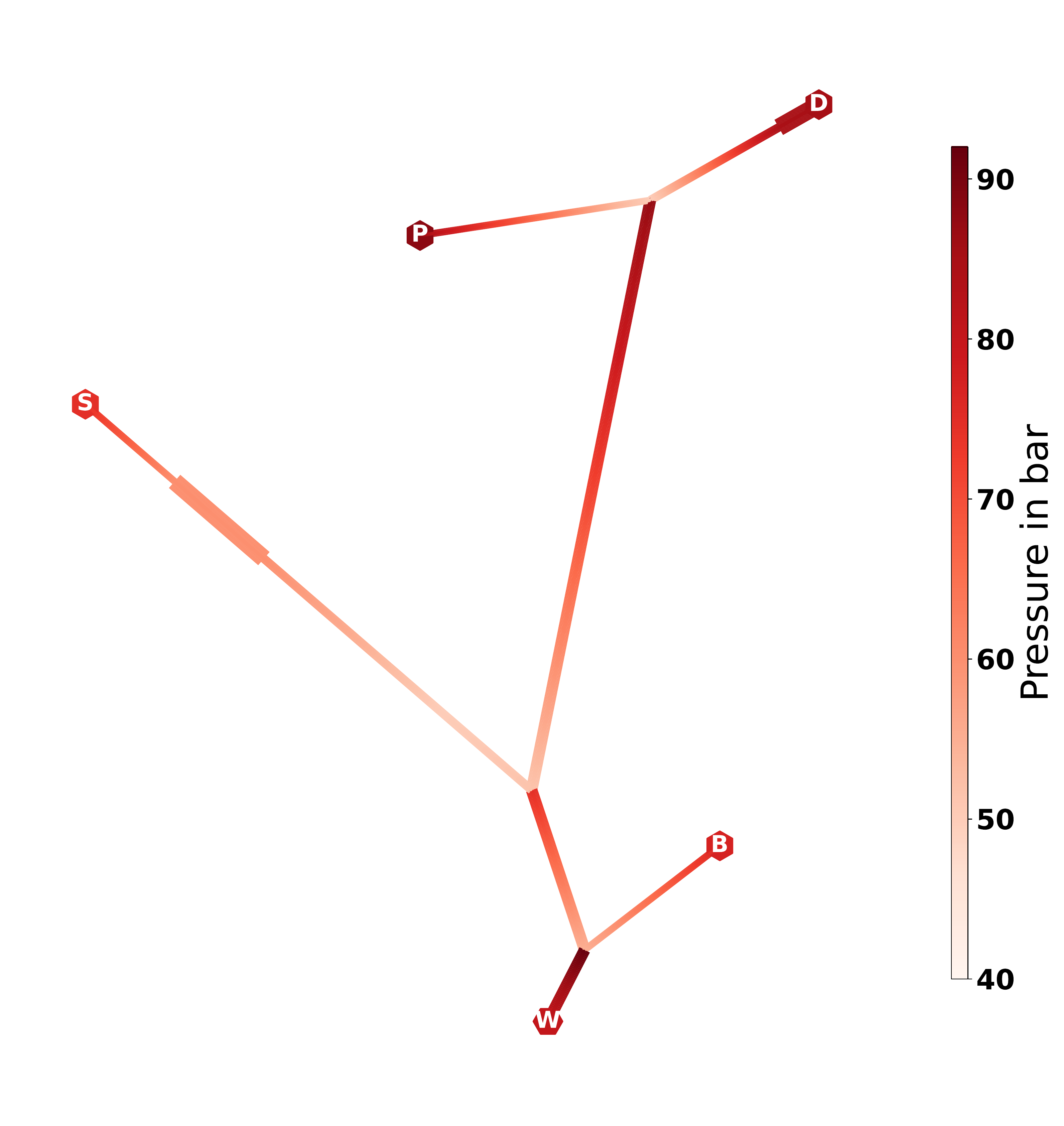}
\caption{Pressure}\label{fig:network_pressure}
\end{subfigure}
\begin{subfigure}[b]{.45\linewidth}
\includegraphics[width=\linewidth]{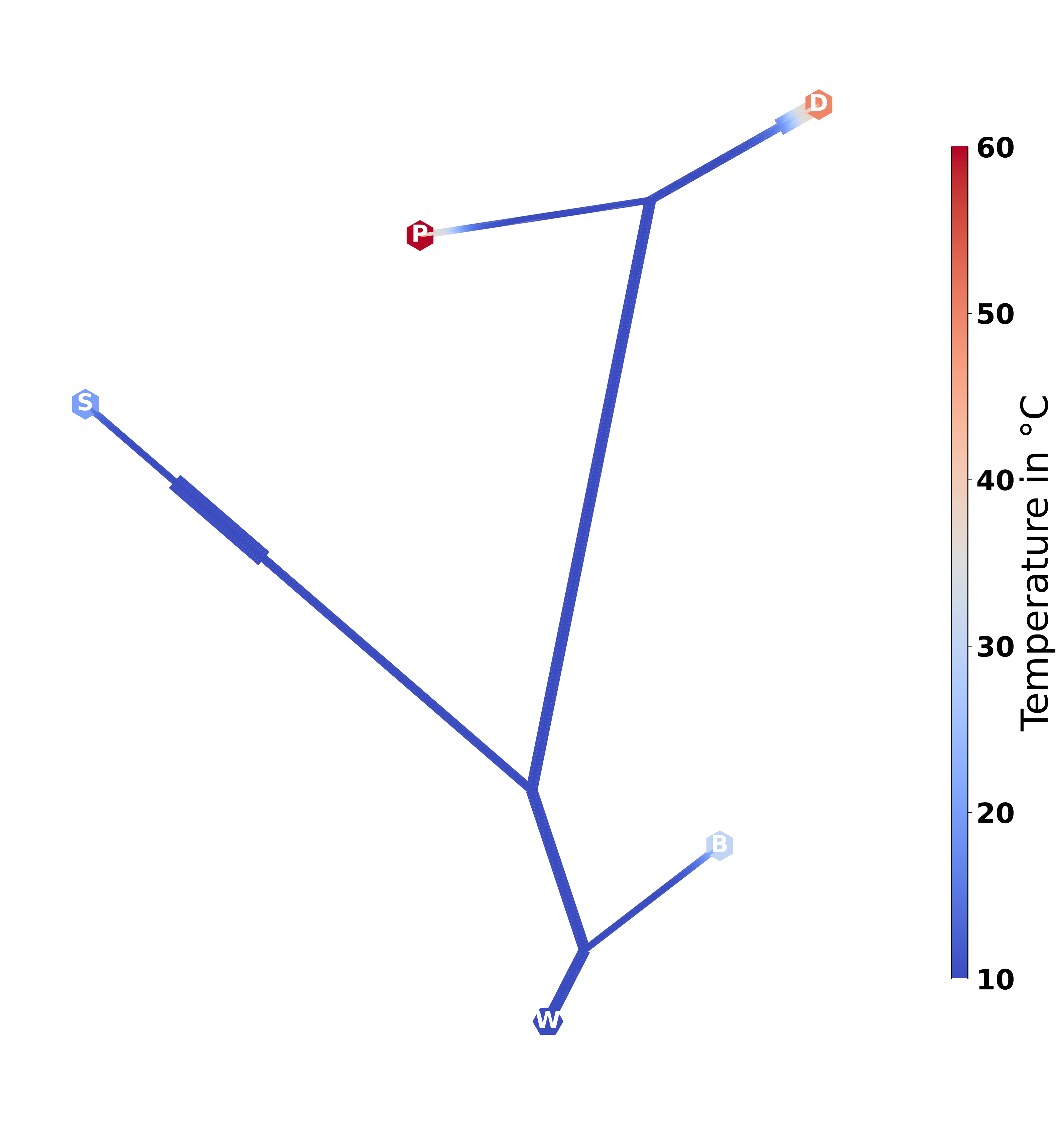}
\caption{Temperature}\label{fig:network_temp}
\end{subfigure}
\caption{Full network result, line thickness proportional to the diameter; Pressure distribution (a) and temperature distribution (b)}
\label{fig:res_network}
\end{figure}

The pressure and temperature distribution are shown in Figure~\ref{fig:network_pressure} and in Figure~\ref{fig:network_temp}, respectively. The inlet pressures at the entry nodes D, P, S, B are at 85.2~\si{bar}, 88.1~\si{bar}, 74.6~\si{bar}, and 77.1~\si{bar}, respectively. The highest pressure values in the network are at 92~\si{bar}, i.e. the pump's outlet pressure. As expected, the \coTwo cools down to the surrounding soil temperature in the first couple of pipes. 

\begin{figure}[htb]
    \centering
    \includegraphics[width=\linewidth]{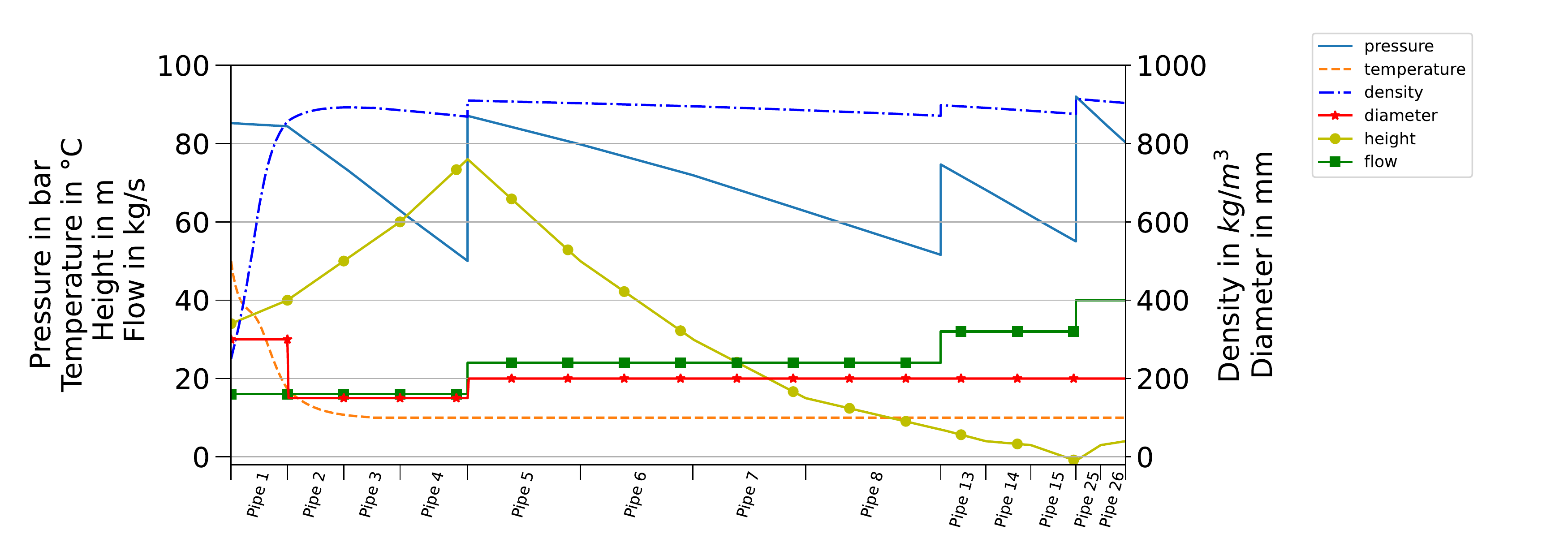}
    \caption{Pressure, temperature, height, flow (left y-axis), density and diameter (right y-axis) profiles for path from entry D to exit W}
    \label{fig:d}
\end{figure}

Figure~\ref{fig:d} shows the profiles of pressure, temperature, height, flow, density, and diameter for the path from entry D to exit node W. The steps in the pressure, density, and flow profiles denote the pumping stations, which are also the junction points at which the different branches of the tree come together. Due to the high temperature difference between the fluid and the surrounding soil, there is a significant temperature drop in the first pipeline. In the middle of Pipe 1, the temperature's graph flattens before declining faster again. This behavior may be caused due to the transition of the fluid from the supercritical into the liquid state (pressure goes below the critical $\pressure_c$ and temperature stays above $\temperature_c$). After Pipe 4, a larger pipeline diameter is needed, as flow from entry node P is added. Despite increasing flow after Pipe 8 and 15, the diameter does not change. Here, the next pumping station or exit is close enough such that the pressure does not drop below the bounds.

\begin{figure}
    \centering
    \includegraphics[width=\linewidth]{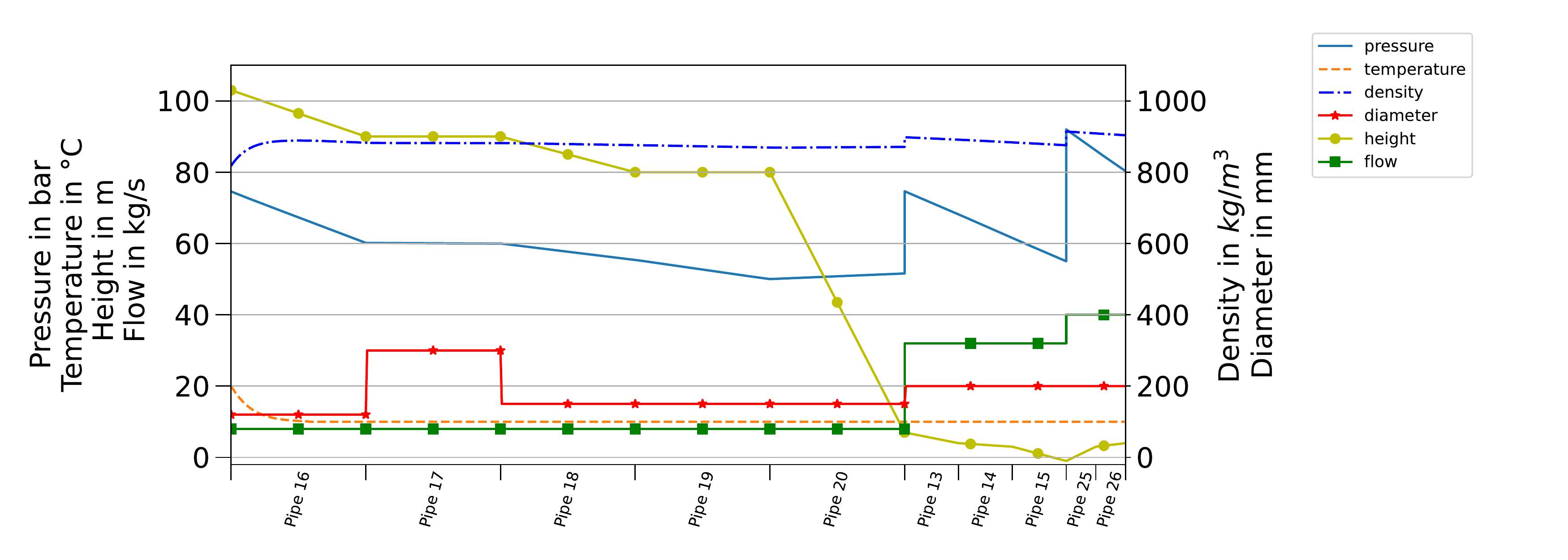}
    \caption{Pressure, temperature, height, flow (left y-axis), density and diameter (right y-axis) profiles for path from entry S to exit W}
    \label{fig:s}
\end{figure}

Figure~\ref{fig:s} shows the profiles of pressure, temperature, height, flow, density and diameter for the path from entry S to the exit node W. Due to the large elevation difference in Pipe 20, the pressure loss because of friction is compensated and there is even a pressure increase. There is almost no pressure loss ($< 0.2$~\si{bar}) on the a priori fixed pipe 17. The statistical profiles of the other entry nodes are shown in \ref{fig:p} and \ref{fig:b} in the appendix.

\subsection{Robustness of the algorithm}\label{sec:robust_analysis}
\begin{figure}
    \centering
    \includegraphics[width=.9\linewidth]{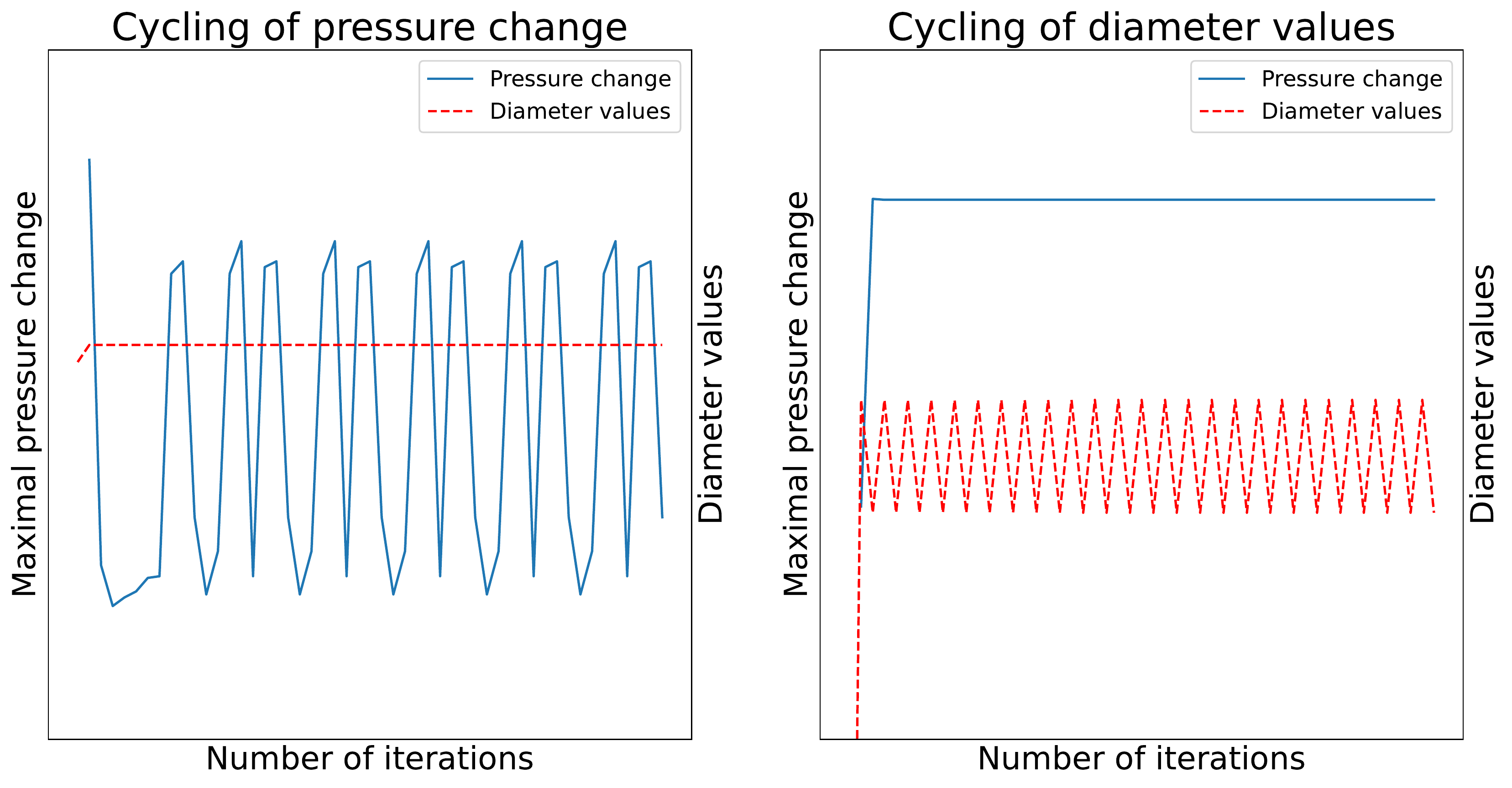}
    \caption{Cycling of pressure change with constant diameter values (left) and cycling of diameter values (right)}
    \label{fig:cycling_example}
\end{figure}

In the following, we test the robustness of our algorithm. In order to do this, we generate a set of artificial tree-shaped network instances of different sizes with various numbers of entry nodes and inflow patterns. Using the data of 120 cities in Germany given in the Travelling Salesman Problem \texttt{gr120.tsp} of the {\sc TSPLIB}95~\cite{reinelt2013}, we randomly choose a number of cities, for which we compute a minimum Steiner Tree based on their distance using the exact solver {\sc SCIP-JACK}~\cite{rehfeldt2021}. Each arc in this tree represents a pipeline. Then, we let the root of the tree be the exit node and all leaves be the entry nodes in the network. For each entry node, we choose an inflow pattern of low (5 \si[per-mode = symbol]{\kg\per\second}, 30\,\% probability), medium (15, 50\,\%), or high (30, 20\,\%) inflow with a low (\SI{20}{\celsius}, 20\,\% probability), medium (\SI{50}{\celsius}, 50\,\%), or high (\SI{80}{\celsius}, 30\,\%) temperature. At each junction node \vertex, i.e., $|\setIncomingArcs{\vertex}|>1$, a pump $\arc^\prime = (\vertex,\vertex^\prime)$ is inserted adding an additional node $\vertex^\prime$ to the network such that the previous outgoing arc of $\vertex$, i.e., $\arc = (\vertex, \anotherVertex)$, is replaced by a pipeline of same length from $\vertex^\prime$ to $\anotherVertex$ and the pump $\arc^\prime$ becomes the new outgoing arc of $\vertex$. We generate 100 instances for 15, 20, 25, and 30 cities each, resulting in 400 instances with the number of entry nodes ranging from three to twelve, the number of pipe segments from 2400 to 6000, and the number of pumps from one to twelve. The computations of this artificial dataset have been performed on an Intel(R) Core(TM) Xeon CPU @ 2.50GHz with 64~GB of RAM. We use as tolerances 
\begin{align*}
    &\frictionTol = 10^{-6}, & \temperatureTol = \SI{0.1}{\kelvin}, && \pressureTol = \SI{1}{\bar}.
\end{align*}
and a maximum number of iteration of 50 as stopping criterion. The computation time is composed out of the time to solve the optimization problems in each iteration and the update step. 306 out of 400 instances are solved to convergence in a mean time of \SI{80}{\second}. For the remaining 94 non-converged instances, we observe two behaviors: One, while the diameters have converged, the change in pressure values decreases until they jump to a higher level only to decrease again, see the left plot in Figure \ref{fig:cycling_example}. Two, the algorithm is cycling between two solutions, one with smaller and one with larger diameters, see the right plot in Figure \ref{fig:cycling_example}. 

\begin{figure}
    \centering
    \includegraphics[width=.7\linewidth]{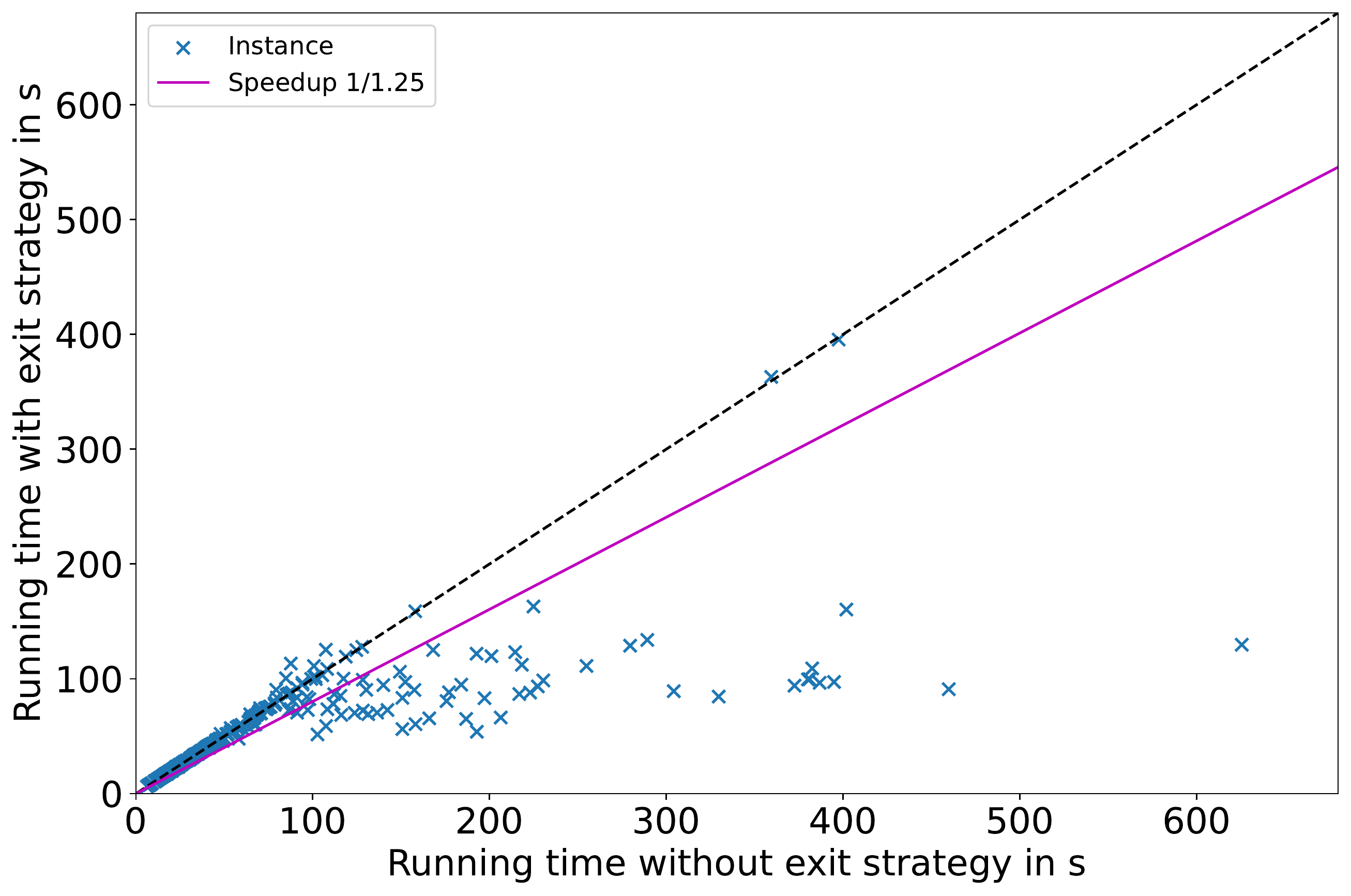}
    \caption{Comparison of running time with and without the proposed exit strategy for solved instances. The black dotted diagonal is the break-even line, the magenta diagonal represents the mean speedup.}
    \label{fig:speed_up}
\end{figure}

In the second case, the solution cycles as a larger pipeline diameter induces a smaller pressure drop which results in an update of the physical quantities, i.e., density, friction factor, et cetera, for which it might be possible to use a smaller diameter and vice versa. Fixing the diameters to the larger values and solving the problem again results in the convergence of the pressure levels. Note that the problem becomes infeasible if the diameters are fixed to the smaller values. For the first case, we propose an exit strategy: If the diameter values do not change for a number of iterations (for our testset, we choose a value of five), we fix the diameters to these values and do not minimize with respect to the diameters, but rather with respect to the inlet pressure of all entries.
This strategy improved the number of converged instances to 372 out of 400. For all instances solved without and with the exit strategy, the final objective values, i.e., choice of diameters, are the same. Figure \ref{fig:speed_up} shows for each instance the computing time with (y-axis) and without (x-axis) the exit strategy. Instances below the dotted diagonal -- the break-even line -- are solved faster with the exit strategy. The magenta diagonal represents the mean speedup of 1.25 for instances solved without and with the exit strategy. In particular, the exit strategy is very effective for instances with larger running times.

\section{Conclusion}\label{sec:conclusion}

In this study, we presented a method to find optimal pipeline diameters in a tree-shaped \coTwo transport network with multiple sources and a single sink taking the physical properties into account. As the behavior of \coTwo is highly sensitive regarding pressure and temperature fluctuations, we propose an iterative algorithm in which the pipeline sizing problem is decoupled from the thermophysical modeling in the network. The problem of finding optimal pipeline diameters from a discrete set of possible diameters is formulated as a MILP for a fixed supply scenario under a given initial temperature distribution. In an update step, the new temperature distribution is determined based on the resulting pipeline diameters and pressure levels in the network. In this step, we incorporate the temperature change due to mixing at junction points and the energy balance of pipelines, including heat exchange with the surrounding soil and the \textit{Joule-Thomson effect}. Afterward, the pipeline sizing problem is solved again with the updated temperatures. After a successful validation on a single pipeline from literature with respect to the pressure drop and temperature profile, the applicability of our approach is shown by solving a real-world \coTwo transport network in North-Western Germany.

Furthermore, the robustness of the approach is validated on a large dataset of artificially created instances. The experiments show that the algorithm may not converge and may cycle. We identify two main reasons for this behavior. One, the pressure does not converge even if the diameters do. And two, the algorithm iterates between two solutions with respect to the choice of diameter. To overcome the first behavior, we propose an exit strategy, changing the objective value to minimal inlet pressure at the entry nodes while we fix the diameters in the network. We have shown the effectiveness of this exit strategy.

There are multiple directions for future extensions to our study. So far, despite the good practical results, theoretical convergence is not proven. Also, even if we have shown optimality with our exit strategy empirically, one could investigate this further. In our study, we focus on pure \coTwo. However, impurities, such as water or nitrogen, have a significant effect, for example, on the fluid's phase envelope and consequently on the pressure bounds, see \cite{nimtz2016}. Additionally, the case of multiple exit nodes could be examined to be able to handle networks featuring multiple places of storage or utilization. Another possibility is to apply our algorithm to other industries, e.g., hydrogen or ammonia transport. Since these are products of demand, e.g., as an energy carrier, one would investigate a single-source multiple-sink network instead of a multiple-source network.

\section{Acknowledgement}
The work for this article has been conducted in the Research Campus MODAL funded by the German Federal Ministry of Education and Research (BMBF) (fund numbers 05M14ZAM, 05M20ZBM).

\bibliographystyle{IEEEtran}
\bibliography{IEEEabrv, CO2network}

\newpage
\appendix
\renewcommand\thefigure{\thesection.\arabic{figure}}
\setcounter{figure}{0}
\renewcommand\thetable{\thesection.\arabic{table}}
\setcounter{table}{0}
\section{Appendix}

\begin{table}[htb]
    \centering
    \begin{tabular}{lccccccc}
     & Pipe 1 & Pipe 2 & Pipe 3 & Pipe 4 & Pipe 5 & Pipe 6 & Pipe 7 \\
    \hline
    Diameter in \si{\milli\meter} & 300& 150& 150& 150& 200& 200& 200\\[0.7ex]
    & Pipe 8 & Pipe 9 & Pipe 10 & Pipe 11 & Pipe 12 & Pipe 13 & Pipe 14 \\
    \hline
    Diameter in \si{\milli\meter} & 200& 120 &120& 120& 120& 200& 200\\[0.7ex]
    & Pipe 15 & Pipe 16 & Pipe 17 &Pipe 18& Pipe 19 & Pipe 20 & Pipe 21 \\ 
    \hline
    Diameter in \si{\milli\meter}& 200& 120& 300& 150 & 150& 150& 120\\[0.7ex]
    & Pipe 22 & Pipe 23 & Pipe 24 & Pipe 25 & Pipe 26 &&\\
    \cline{1-6}
    Diameter in \si{\milli\meter}& 120& 120& 120& 200& 200&&
    \end{tabular}
    \caption{Resulting diameter for each pipeline in the network}
    \label{tab:diameter_sol}
\end{table}

\begin{figure}[htb]
    \centering
    \includegraphics[width=\linewidth]{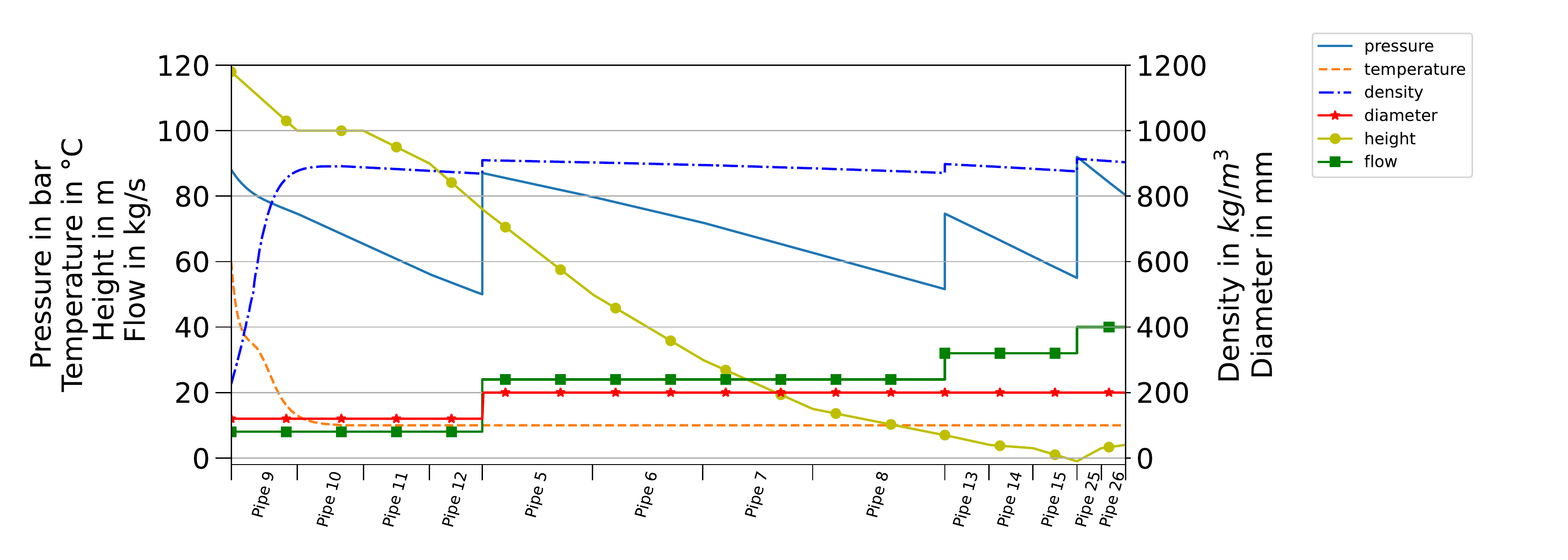}
    \caption{Pressure, temperature, height, flow (left y-axis), density and diameter (right y-axis) profiles for path from entry P to exit W}
    \label{fig:p}
\end{figure}
\begin{figure}[htb]
    \centering
    \includegraphics[width=\linewidth]{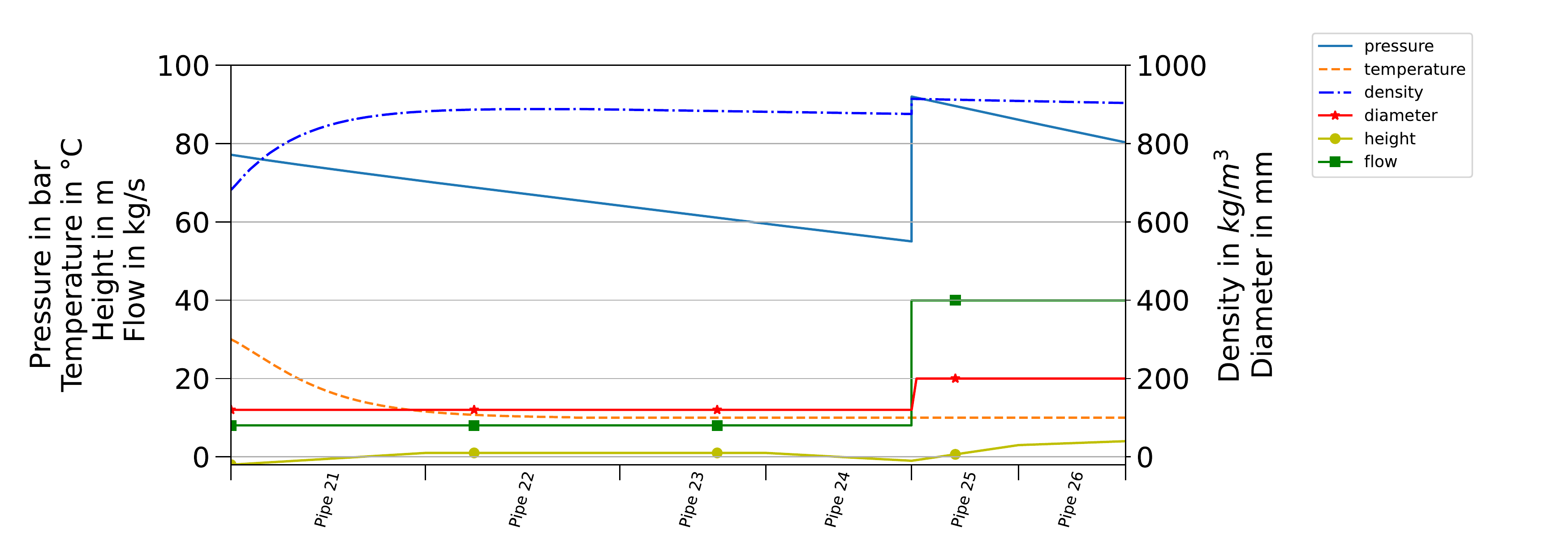}
    \caption{Pressure, temperature, height, flow (left y-axis), density and diameter (right y-axis) profiles for path from entry B to exit W}
    \label{fig:b}
\end{figure}
\end{document}